\documentclass[12pt]{article}
\usepackage{amsmath,amssymb,amsthm}
\def\cc{\mathbb{C}}

\def\nn{{\bf N}}

\def\A{\mathcal A}
\def\B{\mathcal B}
\def\G{\mathcal G}
\def\K{\mathcal K}

\def\F{\mathcal F}

\def\H{\mathcal H}
\def\I{\mathcal I}
\def\J{\mathcal J}
\def\K{\mathcal K}

\def\M{\mathcal M}
\def\N{\mathcal N}
\def\P{\mathcal P}

\def\U{\mathcal U}

\def\amslatex{$\mathcal{A}\kern-.1667em\lower.5ex\hbox{$\mathcal{M}$}\kern-.125em\mathcal{S}$-\LaTeX}

\DeclareMathOperator{\Aut}{Aut}

\def \sp {\hspace{0.3cm}}

\newtheorem{set}{set}[section]
\newtheorem{Corollary}[set]{Corollary}
\newtheorem{Example}[set]{Example}
\newtheorem{Definition}[set]{Definition}

\newtheorem{Lemma}[set]{Lemma}

\newtheorem{Proposition}[set]{Proposition}
\newtheorem{Remark}[set]{Remark}
\newtheorem{Theorem}[set]{Theorem}
\newcommand{\define}{\mathrel{\hbox{$\equiv$\hskip -.90em \lower .47ex \hbox{$\leftharpoondown$}}}}
\newcommand{\enifed}{\mathrel{\hbox{$\equiv$\hskip -.90em \lower .47ex \hbox{$\rightharpoondown$}}}}

\begin{document}
\date{}
\title {Unitarily invariant norms related to factors}
\author{Junsheng Fang \and Don Hadwin}
\maketitle \pagestyle{myheadings}
\begin{abstract} Let $\M$ be a semi-finite factor and
let $\J(\M)$ be the set of operators $T$ in $\M$ such that $T=ETE$
for some finite projection $E$. In this paper we obtain a
representation theorem for unitarily invariant norms on $\J(\M)$
 in terms of Ky Fan norms. As an
application, we prove that the class of unitarily invariant norms
on $\J(\M)$  coincides with the class of symmetric gauge norms on
a classical abelian algebra, which generalizes von Neumann's
classical result~\cite{vN} on unitarily invariant norms on
$M_n(\cc)$. As another application, Ky Fan's dominance
theorem~\cite{Fan} is obtained for  semi-finite factors. Some
classical results in non-commutative $L^p$-theory (e.g.,
non-commutative H$\ddot{\text{o}}$lder's inequality, duality and
reflexivity of non-commutative $L^p$-spaces) are extended to
general unitarily invariant norms related to  semi-finite factors.
We also prove that up to a scale the operator norm is the unique
unitarily invariant norm associated to a type ${\rm III}$ factor.
\end{abstract}



\vskip1.0cm

\section{Introduction}

F.J.~Murray and J.von~Neumann~\cite{M-V1,M-V2,M-V3,vN1,vN2}
introduced and studied certain algebras of Hilbert space
operators. Those algebras are now called ``Von Neumann algebras."
They are strong-operator closed self-adjoint subalgebras of all
bounded linear transformations on a Hilbert space. \emph{Factors}
are von Neumann algebras whose centers consist of scalar multiples
of the identity operator. Every von Neumann algebra is a direct
sum (or ``direct integral'') of factors. Thus factors are the
building blocks for general von Neumann algebras. Murray and von
Neumann~\cite{M-V1} classified factors into type ${\rm I}\sb n,
{\rm I}\sb \infty, {\rm II}\sb 1, {\rm II}\sb \infty, {\rm III}$
factors. Type ${\rm I}\sb n$ and ${\rm I}\sb \infty$ factors are
full matrix algebras: $M_n(\cc)$ and $\B(\H)$. Type ${\rm I}\sb n$
and ${\rm II}\sb 1$ factors are called finite factors. There is a
unique faithful normal tracial state on a finite factor. Factors
except type ${\rm III}$ factors are called semi-finite factors. A
semi-finite factor admits a faithful
normal tracial weight.\\

The unitarily invariant norms on type ${\rm I}\sb n$ factors were
introduced by von Neumann~\cite{vN} for the purpose of metrizing
matrix spaces. Von Neumann, together with his associates,
established that the class of unitarily invariant norms on type
${\rm I}\sb n$ factors coincides with the class of symmetric gauge
norms on $\cc^n$. These norms have now been variously generalized
and utilized in several contexts. For example,
Schatten~\cite{Sc1,Sc2} defined unitarily invariant norms on
two-sided ideals of completely continuous operators in  type ${\rm
I}\sb \infty$ factors; Ky Fan~\cite{Fan} studied Ky Fan norms and
obtained his dominance theorem. The unitarily invariant norms play
a crucial role in the study of function spaces and group
representations (see e.g.~\cite{Ku}) and in obtaining certain
bounds of importance in quantum field theory (see~\cite{Si}). For
historical perspectives and surveys of unitarily invariant norms,
see Schatten~\cite{Sc1,Sc2}, Hewitt and Ross~\cite{H-R}, Gohberg
and
Krein~\cite{G-K} and Simon~\cite{Si}.\\

In~\cite{FHNS}, a structure theorem for unitarily invariant norms
on finite factors is obtained. The main purpose of this paper is
to set up a structure theorem for unitarily invariant norms
related to arbitrary factors, which has a number of
applications.  \\

In this paper, a semi-finite von Neumann algebra $(\M,\tau)$ means
a von Neumann algebra $\M$ with a faithful normal tracial weight
$\tau$, and a Hilbert space $\H$ means the separable
infinite-dimensional complex Hilbert space. If $(\M,\tau)$ is a
finite von Neumann algebra, we assume that $\tau(1)=1$.  If
$\M=\B(\H)$, we assume that $\tau={\rm Tr}$, the classical tracial
weight  on $\B(\H)$.  This paper is organized as the following.\\

In section 2, we define the $s$-numbers of operators in a
semi-finite von Neumann algebra $(\M,\tau)$ from the point of view
of non-increasing rearrangements of functions.\\

In section 3, we study various norms related to a semi-finite von
Neumann algebra $(\M,\tau)$. Let $\J(\M)$ be the set of operators
$T$ in $\M$ such that $T=ETE$ for some finite projection $E$.
 Then $\J(\M)$ is a hereditary
self-adjoint two-sided ideal of $\M$. If $\M$ is a finite von
Neumann algebra, then $\J(\M)=\M$.  If $\M=\B(\H)$, we simply
write $\J(\H)$ instead of $\J(\B(\H))$. Note that $\J(\H)$ is the
set of bounded linear operators $T$ on $\H$ such that both $T$ and
$T^*$ are finite rank operators. A \emph{unitarily invariant} norm
$|\!|\!|\cdot|\!|\!|$ on $\J(\M)$ is a norm on $\J(\M)$ satisfying
$|\!|\!|UTV|\!|\!|=|\!|\!|T|\!|\!|$ for all $T\in \J(\M)$ and
unitary operators $U, V$ in $\M$. For a semi-finite von Neumann
algebra $(\M,\tau)$, let $\Aut(\M,\tau)$ be the set of
$\ast$-automorphisms on $\M$  preserving $\tau$.
 A \emph{symmetric gauge}
 norm $|\!|\!|\cdot|\!|\!|$ on $\J(\M)$ is a  norm on $\J(\M)$
such that $|\!|\!|T|\!|\!|=|\!|\!|\,|T|\,|\!|\!|$ (gauge
invariant) and $|\!|\!|\theta(T)|\!|\!|=|\!|\!|T|\!|\!|$
(symmetric) for all operators $T\in \J(\M)$ and $\theta\in
\Aut(\M,\tau)$.  A norm $|\!|\!|\cdot|\!|\!|$ on $\J(\M)$ is a
\emph{normalized norm} if $|\!|\!|E|\!|\!|=1$ for a projection $E$
in $\M$ such that $\tau(E)=1$.
  We will reserve the notation
$\|\cdot\|$ for the operator norm on a von Neumann algebra.\\

In section 4, we define and study the normalized Ky Fan norms
related to semi-finite von Neumann algebras. To illustrate
difficulties one may encounter in studying of the unitarily
invariant norms related to infinite factors, we point out here one
example. The following result plays a key role in the studying of
unitarily invariant norms on finite factors: if
$|\!|\!|\cdot|\!|\!|$ is a normalized unitarily invariant norm on
a finite factor $(\M,\tau)$, then
\[\|T\|_1\leq |\!|\!|T\!|\!|\leq \|T\|
\] for all $T\in \M$, where $\|T\|_1=\tau(|T|)$ (see Corollary 3.34
of~\cite{FHNS}). However,  the above result is not true for
infinite factors (see Proposition~\ref{P:comparison
proposition}).\\

In section 5, we study the dual norms of symmetric gauge norms on
$\J(\M)$. Let $(\M,\tau)$ be a semi-finite von Neumann algebra and
let $|\!|\!|\cdot|\!|\!|$ be a norm on $\J(\M)$. For $T\in
\J(\M)$, define
\[|\!|\!|T|\!|\!|^\#=\sup\{|\tau(TX)|:\, X\in\J(\M),\, |\!|\!|X|\!|\!|\leq 1\}.
\] In this section, we also compute the dual norms of Ky Fan norms and prove
that  $|\!|\!|\cdot|\!|\!|^{\#\#}=|\!|\!|\cdot|\!|\!|$ under
certain conditions.\\

A representation theorem (Theorem~\ref{T:Theorem C}) for symmetric
gauge norms on $\J(\M)$  is set up in section 6, which is the main
result of this paper. In the rest sections of this paper, we give
a number of
applications of the representation theorem.\\

In section 7, we set up a representation theorem for unitarily
invariant norms related to semi-finite factors and representation
theorems for symmetric gauge norms related to the classical
abelian von Neumann algebras $l^\infty(\mathbf{N})$ and
$L^\infty[0,\infty)$.\\

In section 8, we prove that there is a one-to-one correspondence
between  unitarily invariant norms on $\J(\M)$ for a semi-finite
factor $\M$ and symmetric gauge norms on $\J(\A)$ for a classical
abelian von Neumann algebra $\A$, which generalizes von Neumann's
classical result~\cite{vN} on unitarily invariant norms on type
${\rm I}_n$ factors.  Furthermore, we establish the one to one
correspondence between the dual norms on $\J(\M)$ for a
semi-finite factor $\M$ and the dual norms on $\J(\A)$, which
plays a key role in the studying of duality and reflexivity of the
completion of $\J(\M)$
with respect to  unitarily invariant norms. \\

Ky Fan's dominance theorem for semi-finite factors is proved in
section 9.\\

Let $\M$ be an infinite semi-finite  factor and
$|\!|\!|\cdot|\!|\!|$ be a unitarily invariant norm on $\M$. We
denote by $\overline{\J(\M)_{|\!|\!|\cdot|\!|\!|}}$ the completion
of $\J(\M)$ with respect to $|\!|\!|\cdot|\!|\!|$. Let
$\widetilde{\M}$ be the completion of $\M$ with respect to the
measure topology in the sense of Nelson~\cite{Ne}. In section 10,
we prove that there is an injective map from
$\overline{\J(\M)_{|\!|\!|\cdot|\!|\!|}}$ into $\widetilde{\M}$
that extends the identity map from $\J(\M)$ onto $\J(\M)$. The
reader should compare the proof of Proposition~\ref{P:embedding
into measure topology} and the proof of Proposition 12.3
of~\cite{FHNS} to see the big difference between two cases: type
${\rm II}_1$ case and type ${\rm II}_\infty$ case.
 An element in $\widetilde{\M}$ can be identified with a closed,
densely defined operator affiliated with $\M$ (see~\cite{Ne}). So
generally speaking, an element in
$\overline{\J(\M)_{|\!|\!|\cdot|\!|\!|}}$ should be treated as an
unbounded operator. In section 10, we also extend the
non-commutative H$\ddot{\text{o}}$lder's inequality to the general
unitarily invariant norms. \\

 In section 11, we  study the duality and reflexivity of
  $\overline{\J(\M)_{|\!|\!|\cdot|\!|\!|}}$
 for infinite semi-finite factors. We point out another example of difference
 between two cases: finite factors and infinite semi-finite
 factors.
  Let $\N$ be a type ${\rm II}\sb
1$ factor with a faithful normal state $\tau_\N$. Let
$|\!|\!|\cdot|\!|\!|$ be
 a unitarily invariant norm on $\N$ and let
$|\!|\!|\cdot|\!|\!|^\#$ be the dual unitarily invariant norm on
$\N$. The following
 result is proved in~\cite{FHNS}: $\overline{{\N}_{|\!|\!|\cdot|\!|\!|^\#}}$ is the dual
space of $\overline{{\N}_{|\!|\!|\cdot|\!|\!|}}$ if and only if
$|\!|\!|\cdot|\!|\!|$ is a continuous norm on $\N$, i.e.,
$\lim_{\tau(E)\rightarrow 0+}|\!|\!|E|\!|\!|=0$. A key step to
proving the above result is based on the following fact: if
$|\!|\!|\cdot|\!|\!|$ is a continuous unitarily invariant norm on
$\N$ and $\phi\in \overline{\N_{|\!|\!|\cdot|\!|\!|}}^\#$, then
the restriction of $\phi$ to $\N$ is an ultraweakly continuous
linear functional, i.e., $\phi$ is in the predual space of $\N$.
However, it is easy to see that the similar result is not true for
infinite semi-finite factors $\M$, e.g., $\M=\B(\H)$.

\vskip 0.5cm
 \noindent \emph{Acknowledgements:}\, In the two
semesters from 2006-2007, the first author learned a lot of
mathematics from Professor Eric Nordgren. The first author expresses
his deep gratitude to Eric for his teaching,  encouragements and
many valuable suggestions.
\section{Preliminaries}
\subsection{Nonincreasing rearrangements of functions}
Throughout this paper, we denote by $m$ the Lebesgue measure on
$[0,\infty)$. In the following, a measurable function and a
measurable set mean a Lesbesgue measurable function and a Lebesgue
measurable set, respectively.  Let  $f(x)$ be a real measurable
function on $[0,\infty)$. The \emph{nonincreasing rearrangement
function}, $f^*(x)$, of $f(x)$ is defined by
\begin{equation}\label{E:nonincreasing}
f^*(x)=\sup\{y:\, m(\{f>y\})>x\}, \quad 0\leq x<\infty.
\end{equation}
We summarize some well-known properties of $f^*(x)$ in the following
proposition.

\begin{Proposition}\label{P:nonincreasing} Let $f(x), g(x), f_1(x), f_2(x), \cdots$ be
real measurable functions on $[0,\infty)$, $c$ be a real number.
Then we have the following:
 \begin{enumerate}
\item  $f^*(x)$ is a nonincreasing, right-continuous function on
$[0,\infty)$ such that $f^*(0)=essup f(x)$;
\item  $(f+c)^*=f^*+c$;
\item $(cf)^*=cf^*$ if $c\geq 0$;
\item if $f(x)$ is a simple function, then $f^*(x)$ is also a
simple function;
\item if $f(x)\leq g(x)$ for almost all $x$, then $f^*(x)\leq
g^*(x)$ everywhere;
\item  $\|f^*(x)-g^*(x)\|_\infty\leq \|f(x)-g(x)\|_\infty$;
\item if $\lim_{n\rightarrow\infty}f_n(x)=f(x)$ uniformly, then
$\lim_{n\rightarrow\infty} f_n^*(x)=f^*(x)$ uniformly;
\item if $f_n(x)$ converges to $f(x)$ in measure, then
$\liminf_{n\rightarrow\infty}f_n^*(x)\geq f^*(x)$ for every $x\in
[0,\infty)$;
\item if $f_n(x)$ converges to $f(x)$ in measure, then $\limsup_{n\rightarrow\infty}f_n^*(x)\leq f^*(x)$ for every $x\in
[0,\infty)$ such that $f^*$ is continuous at $x$;
\item $f(x)$ and $f^*(x)$ are equi-measurable, i.e, for any real
number $y$, $m(\{f>y\})=m(\{f^*>y\})$;
\item $f^*=g^*$ if and only $f(x)$ and $g(x)$ are equi-measurable;
\item if $f(x)$ and $g(x)$ are bounded  functions and
$\int_0^\infty f^n(x)dx=\int_0^\infty g^n(x)dx$ for all
$n=0,1,2,\cdots,$ then $f^*(x)=g^*(x)$;
\item $\int_0^\infty f(x)dx=\int_0^\infty f^*(x)dx$ when either  integral is
well-defined.
\end{enumerate}
\end{Proposition}

\subsection{$s$-numbers of operators in type ${\rm II}\sb \infty$ factors}
In~\cite{F-K}, Fack and Kosaki give a rather complete exposition
of generalized $s$-numbers  of $\tau-$measurable operators
affiliated with  semi-finite von Neumann algebras. For the
reader's convenience and our purpose, we provide sufficient
details on $s$-numbers of bounded operators in semi-finite von
Neumann algebras in this section. We will define $s$-numbers of
bounded operators in semi-finite von Neumann algebras from the
point of view of non-increasing rearrangements of functions. The
following lemma is well-known.

\begin{Lemma}\label{L:isomorphism}
 Let $(\A,\tau)$ be a separable \emph{(}i.e., with separable
 predual\emph{)}
 diffuse abelian von Neumann algebra  with a faithful normal tracial weight $\tau$ on $\A$
 such that $\tau(1)=\infty$. Then there is
a $*-$isomorphism $\alpha$ from $(\A,\tau)$ onto
$(L^\infty[0,\infty),\int_0^\infty dx)$ such that
$\tau=\int_0^\infty dx\cdot \alpha$.
\end{Lemma}

Let $\M$ be a type ${\rm II}\sb \infty$ factor and let $\tau$ be a
faithful normal tracial weight on $\M$. For $T\in \M$, there is a
separable diffuse abelian von Neumann subalgebra $\A$ of $\M$
containing $|T|$. By Lemma~\ref{L:isomorphism}, there is a
$\ast$-isomorphism $\alpha$ from $(\A,\tau)$ onto
$(L^\infty[0,\infty), \int_0^\infty dx)$ such that
$\tau=\int_0^\infty dx\cdot\alpha$. Let $f(x)=\alpha(|T|)$ and let
$f^*(x)$ be the non-increasing rearrangement of $f(x)$
(see~$(\ref{E:nonincreasing})$). Then \emph{the $s$- numbers of
$T$}, $\mu_s(T)$, are defined as
\[\mu_s(T)=f^*(s),\,\, 0\leq s<\infty.
\]

\begin{Lemma}\label{L:s-number}$\mu_s(T)$ does not depend on $\A$ and $\alpha$.
\end{Lemma}
\begin{proof} Let $\A_1$ be another separable diffuse abelian von
Neumann subalgebra of $\M$ containing $|T|$ and suppose $\beta$ is
a $\ast$-isomorphism from $\A_1$ onto $L^\infty[0,\infty)$ such
that $\tau=\int_0^\infty dx\cdot\beta$. Let $g(x)=\beta(|T|)$. For
every number $n=0,1,2,\cdots$, $\int_0^\infty
f^n(x)dx=\tau(|T|^n)=\int_0^\infty g^n(x)dx$. Since both $f(x)$
and $g(x)$ are bounded positive functions, by 12 of
Proposition~\ref{P:nonincreasing}, $f^*(x)=g^*(x)$ for all $x\in
[0,\infty)$.
\end{proof}
\begin{Corollary}\label{C:s-number} For $T\in \M$ and $p\geq 0$, $\tau(|T|^p)=\int_0^\infty \mu_s(T)^pds$.
\end{Corollary}

\begin{Lemma}\label{L:projections} Let $E,F$ be two projections in
$\M$. If $\tau(E^\perp)<\tau(F^\perp)<\infty$, then $\tau(E\wedge
F^\perp)>0$.
\end{Lemma}
\begin{proof} By Proposition~2.5.14  of~\cite{F-K} (page 119, vol
1), $R(F^\perp E^\perp)=F^\perp-E\wedge F^\perp$, where $R(F^\perp
E^\perp)$ is the range projection of $F^\perp E^\perp$. Therefore,
\[\tau(E\wedge F^\perp)=\tau(F^\perp)-\tau(R(F^\perp E^\perp))\geq
\tau(F^\perp)-\tau(E^\perp)>0.\]

\end{proof}

Let $\P(\M)$ be the set of projections in $\M$. The following
lemma says that above definition of $s$-numbers coincides with the
definition of $s$-numbers given by Fack and Kosaki.

\begin{Lemma}\label{L:s-number2} For $0\leq s<\infty$,
\[\mu_s(T)=\inf\{\|TE\|: \,\, E\in\P(\M),\, \tau(E^\perp)=s\}.
\]
\end{Lemma}
\begin{proof} By the polar decomposition and the definition of
$\mu_s(T)$, we may assume that $T$ is positive. Let $\A$ be a
separable diffuse abelian von Neumann subalgebra of $\M$
containing $T$ and let $\alpha$ be a $\ast$-isomorphism from $\A$
onto $L^\infty[0,\infty)$ such that $\tau=\int_0^\infty\, dx\cdot
\alpha$. Let $f(x)=\alpha(T)$ and let $f^*(x)$ be the
nonincreasing rearrangement of $f(x)$. Then $\mu_s(T)=f^*(s)$. By
the definition of $f^*$, \[m(\{f^*>\mu_s(T)\})=
\lim_{n\rightarrow\infty}m\left(\left\{f^*>\mu_s(T)+\frac{1}{n}\right\}\right)\leq
s\] and
\[m(\{f^*\geq \mu_s(T)\})\geq \lim_{n\rightarrow\infty} m\left(\left\{f^*>
\mu_s(T)-\frac{1}{n}\right\}\right)\geq s.\] Since $f^*$ and $f$ are
equi-measurable, $m(\{f>\mu_s(T)\})\leq s$ and $m(\{f\geq
\mu_s(T)\})\geq s$. Therefore, there is a measurable set $A$ of
$[0,\infty)$, $\{f>\mu_s(T)\}\subset [0,\infty)\setminus A\subset
\{f\geq \mu_s(T)\}$, such that $m([0,\infty)\setminus A)=s$ and
$\|f(x)\chi_A(x)\|_\infty=\mu_s(T)$ and
$\|f(x)\chi_B(x)\|_\infty\geq \mu_s(T)$ for every $B\subset
[0,\infty)\setminus A$ such that $m(B)>0$. Let
$F=\alpha^{-1}(\chi_A)$. Then $\tau(F^\perp)=s$,
$\|TF\|=\|\alpha^{-1}(f\chi_A)\|_\infty=\mu_s(T)$ and $\|TF'\|\geq
\mu_s(T)$ for any nonzero subprojection $F'$ of $F^\perp$. This
proves that $\mu_s(T)\geq \inf\{\|TE\|: \,\, E\in\P(\M),\,
\tau(E^\perp)=s\}$. Similarly, for any $\epsilon>0$, there is a
projection $F_\epsilon\in \M$ such that
$\tau(F_\epsilon^\perp)=s+\epsilon$,
$\|TF_\epsilon\|=\mu_{s+\epsilon}(T)$ and $\|TF'\|\geq
\mu_{s+\epsilon}(T)$ for any nonzero subprojection $F'$ of
$F_\epsilon^\perp$.  Suppose $E\in \M$ is a projection such that
$\tau(E^\perp)=s$. By Lemma~\ref{L:projections}, $\tau(E\wedge
F_\epsilon^\perp)>0$. Hence, $\|TE\|\geq \|T(E\wedge
F_\epsilon^\perp)\|\geq \mu_{s+\epsilon}(T)$. This proves that
$\inf\{\|TE\|: \,\, E\in\P(\M),\, \tau(E^\perp)=s\}\geq
\mu_{s+\epsilon}(T)$. Since $\mu_s(T)$ is right-continuous,
$\mu_s(T)\leq \inf\{\|TE\|: \,\, E\in\P(\M),\, \tau(E^\perp)=s\}.$
\end{proof}

\begin{Corollary}\label{C:property of s-number} Let $S,T\in \M$. Then $\mu_s(ST)\leq
\|S\|\mu_s(T)$ for $s\in [0,\infty)$.
\end{Corollary}

We refer to~\cite{Fa,F-K} for other interesting properties of
$s$-numbers for operators in type ${\rm II}\sb \infty$ factors.
\subsection{$s$-numbers of operators in  semi-finite von Neumann
algebras}

 An \emph{embedding}
of a semi-finite von Neumann algebra $(\M,\tau)$ into another
semi-finite von Neumann algebra $(\M_1,\tau_1)$ means a
$\ast$-isomorphism $\alpha$ from $\M$ to $\M_1$ such that
$\tau=\tau_1\cdot\alpha$. Every  semi-finite von Neumann algebra can
be embedded into a type  ${\rm II}\sb \infty$ factor.

\begin{Definition}\label{D:s-number}\emph{ Let $(\M,\tau)$ be a semi-finite von
Neumann algebra and $T\in \M$. If $\alpha$ is an embedding of
$(\M,\tau)$ into a type ${\rm II}\sb \infty$ factor
$(\M_1,\tau_1)$, then \emph{the $s$-numbers of $T$} are defined as
\[\mu_s(T)=\mu_s(\alpha(T)),\quad 0\leq s<\infty.
\]}
\end{Definition}

Similar to the proof of Lemma~\ref{L:s-number}, we can see that
$\mu_s(T)$ is well defined, i.e., does not depend on the choice of
$\alpha$ and $\M_1$.\\

 Let $T\in (\B(\H),{\rm Tr})$ be a finite rank operator,
where $\H$ is the separable infinite dimensional complex Hilbert
space and ${\rm Tr}$ is the classical tracial weight on $\B(\H)$.
Then $|T|$ is unitarily equivalent to a diagonal operator  with
diagonal elements $s_1(T)\geq s_2(T)\geq \cdots\geq 0$. In the
classical operator theory~\cite{G-K}, $s_1(T), s_2(T), \cdots$ are
also called $s$-numbers of $T$. It is easy to see that the relation
between $\mu_s(T)$ and $s_1(T),s_2(T),\cdots$ is the following
\begin{equation}\label{E:s-number}
\mu_s(T)=s_1(T)\chi_{[0,1)}(s)+s_2(T)\chi_{[1,2)}(s)+\cdots.
\end{equation}
Since no confusion will arise, we will use both $s$-numbers for a
finite rank operator in $(\B(\H),{\rm Tr})$. We refer to~\cite{Bh,
G-K} for properties of $s$-numbers for finite rank operators in
$(\B(\H),{\rm Tr})$.

We end this section by the following definition.
\begin{Definition}\label{D:equimeasurable}\emph{ Two positive operators $S,T$ in a semi-finite von Neumann
algebra $(\M,\tau)$ are \emph{equi-measurable} if
$\mu_s(S)=\mu_s(T)$ for $0\leq s<\infty$.}
\end{Definition}

By 12 of Proposition~\ref{P:nonincreasing} and
Corollary~\ref{C:s-number}, positive operators $S$ and $T$ in a
semi-finite von Neumann algebra $(\M,\tau)$ are equi-measurable if
and only if $\tau(S^n)=\tau(T^n)$ for all $n=0,1,2,\cdots$.

\section{Semi-norms on $\J(\M)$}

In this section, $(\M,\tau)$ is a semi-finite von Neumann algebra
with a faithful normal tracial weight $\tau$. Recall that $\J(\M)$
is the set of operators $T$ in $\M$ such that $T=ETE$ for some
finite projection $E$.  If $\M=\B(\H)$, we simply write $\J(\H)$
instead of $\J(\B(\H))$. Note that $\J(\H)$ is the set of bounded
linear operators $T$ on $\H$ such that both $T$ and $T^*$ are
finite rank operators. \\

\subsection{Gauge invariant semi-norms on $\J(\M)$}

\begin{Definition}\label{D:gauge invariant norm on J(M)}
\emph{ Let $(\M,\tau)$ be a  semi-finite  von Neumann algebra.
 A semi-norm $|\!|\!|\cdot|\!|\!|$ on $\J(\M)$ is \emph{gauge invariant}
  if $|\!|\!|T|\!|\!|=|\!|\!|\,|T|\,|\!|\!|$ for
all $T\in \J(\M)$. A semi-norm $|\!|\!|\cdot|\!|\!|$ on $\J(\M)$ is
called \emph{left unitarily invariant} if for all unitary operators
$U$ in $\M$ and all $T$ in $\J(\M)$,
$|\!|\!|UT|\!|\!|=|\!|\!|T|\!|\!|.$}
\end{Definition}

\begin{Lemma}\label{L:submultiplicative} Let $(\M,\tau)$ be a semi-finite von Neumann algebra
and let $|\!|\!|\cdot|\!|\!|$ be a left unitarily invariant
semi-norm on $\J(\M)$. If $T\in \J(\M)$ and $A\in \M$, then $AT\in
\J(\M)$ and $|\!|\!|AT|\!|\!|\leq \|A\|\cdot|\!|\!|T|\!|\!|$.
\end{Lemma}
\begin{proof} It is easy to see that $AT\in \J(\M)$. We need to prove that if
$\|A\|<1$, then $|\!|\!|AT|\!|\!|\leq |\!|\!|T|\!|\!|$. Since
$\|A\|<1$, there are unitary operators $U_1,\cdots,U_k$ such that
$A=\frac{U_1+\cdots+U_k}{k}$ (see~\cite{K-P, R-D}). Since
$|\!|\!|\cdot|\!|\!|$ is a left unitarily invariant semi-norm on
$\J(\M)$, \[|\!|\!|AT|\!|\!|=
|\!|\!|\frac{U_1T+\cdots+U_kT}{k}|\!|\!|\leq
\frac{|\!|\!|U_1T|\!|\!|+\cdots+|\!|\!|U_kT|\!|\!|}{k}\leq
|\!|\!|T|\!|\!|.\]
\end{proof}

\begin{Lemma}\label{L:gauge=left unitarily invariant}Let $(\M,\tau)$
 be a semi-finite von Neumann algebra and  let $|\!|\!|\cdot|\!|\!|$ be a semi-norm
on $\J(\M)$. Then $|\!|\!|\cdot|\!|\!|$ is gauge invariant if and
only if $|\!|\!|\cdot|\!|\!|$ is left unitarily invariant.
\end{Lemma}
\begin{proof} Note that
$|UT|=|T|$ for all $T\in \J(\M)$ and unitary operators $U$ in $\M$.
If $|\!|\!|\cdot|\!|\!|$ is gauge invariant then
$|\!|\!|\cdot|\!|\!|$ is left unitarily invariant. Conversely,
suppose $|\!|\!|\cdot|\!|\!|$ is left unitarily invariant. By
Lemma~\ref{L:submultiplicative},
$|\!|\!|T|\!|\!|=|\!|\!|V|T|\,\,|\!|\!|\leq |\!|\!|\,\,
|T|\,\,|\!|\!|$ and $|\!|\!|\,\,
|T|\,\,|\!|\!|=|\!|\!|V^*T|\!|\!|\leq |\!|\!|T|\!|\!|.$ Hence,
$|\!|\!|\cdot|\!|\!|$ is gauge invariant.
\end{proof}

\begin{Corollary}\label{C:S<T}Let $(\M,\tau)$ be a semi-finite  von Neumann algebra and let $|\!|\!|\cdot|\!|\!|$ be
a gauge invariant semi-norm on $\J(\M)$. If $T\in \J(\M)$ and $0\leq
S\leq T$, then $S\in \J(\M)$ and $|\!|\!|S|\!|\!|\leq
|\!|\!|T|\!|\!|$.
\end{Corollary}
\begin{proof} Since $0\leq S\leq T$, there is an operator $A\in \M$
such that $S=AT$ and $\|A\|\leq 1$. By
Lemma~\ref{L:submultiplicative} and Lemma~\ref{L:gauge=left
unitarily invariant}, $S\in \J(\M)$ and
$|\!|\!|S|\!|\!|=|\!|\!|AT|\!|\!|\leq
\|A\|\cdot|\!|\!|T|\!|\!|\leq |\!|\!|T|\!|\!|$.
\end{proof}

\subsection{Unitarily invariant semi-norms on $\J(\M)$}

\begin{Definition}\label{D:unitarily invariant norm on J(M)}
\emph{ Let $(\M,\tau)$ be a  semi-finite  von Neumann algebra.
 A semi-norm $|\!|\!|\cdot|\!|\!|$ on $\J(\M)$ is \emph{unitarily invariant}
  if $|\!|\!|UTV|\!|\!|=|\!|\!|T|\!|\!|$ for
all $T\in \J(\M)$ and unitary operators $U,V\in \M$.}
\end{Definition}

\begin{Proposition}\label{P:unitarily invariant norms on J(M)} Let $|\!|\!|\cdot|\!|\!|$ be a semi-norm on $\J(\M)$. Then the
following statements are equivalent:
\begin{enumerate}
\item $|\!|\!|\cdot|\!|\!|$ is unitarily invariant;
\item $|\!|\!|\cdot|\!|\!|$ is gauge invariant and unitarily conjugate
invariant, i.e.,  $|\!|\!|UTU^*|\!|\!|=|\!|\!|T|\!|\!|$ for all
$T\in\J(\M)$ and unitary operators $U\in \M$;
\item $|\!|\!|\cdot|\!|\!|$ is
gauge invariant and $|\!|\!|T|\!|\!|=|\!|\!|T^*|\!|\!|$ for all
$T\in \J(\M)$;
\item for all operators
$A,B\in \M$ and $T\in \J(\M)$, $|\!|\!|ATB|\!|\!|\leq \|A\|\cdot
|\!|\!|T|\!|\!|\cdot \|B\|$.
\end{enumerate}
\end{Proposition}
\begin{proof} $``1\Rightarrow 4"$ is similar to the proof of Lemma~\ref{L:submultiplicative}.

$``4\Rightarrow 3"$. Let $T=V|T|$. Then $T^*=|T|V^*$. By 4 and
simple arguments, $|\!|\!|T|\!|\!|=|\!|\!|T^*|\!|\!|$.

$``3\Rightarrow 2"$. By Lemma~\ref{L:gauge=left unitarily invariant}
and 3,
$|\!|\!|UTU^*|\!|\!|=|\!|\!|TU^*|\!|\!|=|\!|\!|UT^*|\!|\!|=|\!|\!|T^*|\!|\!|=
|\!|\!|T|\!|\!|$.

$``2\Rightarrow 1"$.
 Suppose $|\!|\!|\cdot|\!|\!|$ is
gauge invariant and unitarily conjugate invariant. Let $U,V\in \M$
be unitary operators and $T\in \J(\M)$. By Lemma~\ref{L:gauge=left
unitarily invariant},
$|\!|\!|UTV|\!|\!|=|\!|\!|V^*VUTV|\!|\!|=|\!|\!|VUT|\!|\!|=|\!|\!|T|\!|\!|$.

\end{proof}

\begin{Corollary}\label{C:equivalent projections in J(M)} Let $|\!|\!|\cdot|\!|\!|$
be a unitarily invariant semi-norm on $\J(\M)$ and let $E, F$ be
two equivalent projections in $\J(\M)$. Then
$|\!|\!|E|\!|\!|=|\!|\!|F|\!|\!|$.
\end{Corollary}

\subsection{Symmetric gauge semi-norms on $\J(\M)$}
\begin{Definition}\label{D:Symmetric gauge norms}
 \emph{Let $(\M,\tau)$ be a semi-finite von Neumann algebra  and let $\Aut(\M, \tau)$ be  the set of
$\ast$-automorphisms  on $\M$ preserving $\tau$.
 A semi-norm
 $|\!|\!|\cdot|\!|\!|$  on $\J(\M)$ is called \emph{symmetric} (with respect to $\tau$)  if
\[|\!|\!|\theta(T)|\!|\!|=|\!|\!|T|\!|\!|,\,\,\,\, \forall\, T\in \J(\M),\, \theta\in
\Aut(\M,\tau);
\]  a semi-norm
 $|\!|\!|\cdot|\!|\!|$  on $\J(\M)$ is called a \emph{symmetric gauge semi-norm}  if it
is both symmetric and gauge invariant on $\J(\M)$.}
\end{Definition}

\begin{Example}\emph{The abelian von Neumann algebra
$\mathbb{C}^{n}$ is a finite von Neumann algebra with a classical
tracial state $\tau((x_1,\cdots,x_n))=\frac{x_1+\cdots+x_n}{n}$.
In this case, $\J(\mathbb{C}^{n})=\mathbb{C}^{n}$. A norm
$|\!|\!|\cdot|\!|\!|$ on $\mathbb{C}^{n}$ is a symmetric gauge
norm if and only if for every $(x_1,\cdots, x_n)\in
\mathbb{C}^{n}$,
\begin{enumerate}
\item $|\!|\!|(x_1,\cdots, x_2)|\!|\!|=|\!|\!|(|x_1|,\cdots,
|x_n|)|\!|\!|$ and \item $|\!|\!|(x_1,\cdots,
x_n)|\!|\!|=|\!|\!|(x_{\pi(1)},\cdots, x_{\pi(n)})|\!|\!|$ for
every permutation $\pi$ of $\{1,\cdots,n\}$.
\end{enumerate} }
\end{Example}

\begin{Example}\emph{The abelian von Neumann algebra
$l^\infty(\mathbf{N})$ is a semi-finite von Neumann algebra with a
classical tracial weight $\tau((x_1,x_2,\cdots))=x_1+x_2+\cdots$.
It is easy to see that $\J(l^\infty(\mathbf{N}))=c_{00}$ consists
of $(x_1,x_2,\cdots)$ with $x_n=0$ except for finitely many $n$. A
norm $|\!|\!|\cdot|\!|\!|$ on $\J(l^\infty(\mathbf{N}))$ is a
symmetric gauge norm if and only if for every $(x_1,x_2,\cdots)\in
c_{00}$,
\begin{enumerate}
\item
$|\!|\!|(x_1,x_2,\cdots)|\!|\!|=|\!|\!|(|x_1|,|x_2|,\cdots)|\!|\!|$
and \item
$|\!|\!|(x_1,x_2,\cdots)|\!|\!|=|\!|\!|(x_{\pi(1)},x_{\pi(2)},\cdots)|\!|\!|$
for every permutation $\pi$ of $\mathbf{N}$.
\end{enumerate} }
\end{Example}

\begin{Example}\emph{The abelian von Neumann algebra
$L^\infty[0,1]$ is a finite von Neumann algebra with a classical
tracial state $\tau=\int_0^1\,dx$. In this case
$\J(L^\infty[0,1])=L^\infty[0,1]$.  A norm $|\!|\!|\cdot|\!|\!|$
on $L^\infty[0,1]$ is a symmetric gauge norm  if and only if for
every $f(x)\in L^\infty[0,1]$,
\begin{enumerate}
\item $|\!|\!|f(x)|\!|\!|=|\!|\!|\,|f(x)|\,|\!|\!|$ and \item
$|\!|\!|f(x)|\!|\!|=|\!|\!|f(\phi(x))|\!|\!|$ for every invertible
measure preserving map $\phi$ of $[0,1]$.
\end{enumerate} }
\end{Example}

\begin{Example}\emph{The abelian von Neumann algebra
$L^\infty[0,\infty)$ is a semi-finite von Neumann algebra with a
classical tracial weight $\tau=\int_0^\infty\,dx$.  A norm
$|\!|\!|\cdot|\!|\!|$ on $\J(L^\infty[0,\infty))$ is a symmetric
gauge norm  if and only if for every
$f(x)\in\J(L^\infty[0,\infty))$,
\begin{enumerate}
\item $|\!|\!|f(x)|\!|\!|=|\!|\!|\,|f(x)|\,|\!|\!|$ and \item
$|\!|\!|f(x)|\!|\!|=|\!|\!|f(\phi(x))|\!|\!|$ for every invertible
measure preserving map $\phi$ of $[0,\infty)$.
\end{enumerate} }
\end{Example}

\begin{Lemma}\label{C:symmetric implies unitarily invariant on J(M)} Let
$(\M,\tau)$ be a semi-finite von Neumann algebra and let
$|\!|\!|\cdot|\!|\!|$ be a symmetric gauge semi-norm on $\J(\M)$.
Then $|\!|\!|\cdot|\!|\!|$ is  a  unitarily invariant semi-norm on
$\J(\M)$.
\end{Lemma}

\subsection{Symmetric gauge norms on $(\M_E,\tau_E)$}
 In this paper we are
interested symmetric gauge semi-norms on $\J(\M)$, where $(\M,\tau)$
is one of the following semi-finite von Neumann algebras:
\begin{itemize}
\item $\M=\B(\H)$ and $\tau={\rm Tr}$ on $\B(\H)$, where $\H$ is the separable infinite dimensional complex Hilbert space;
\item $\M=l^\infty(\nn)$ and
$\tau((x_1,x_2,\cdots))=x_1+x_2+\cdots$;
\item $\M$ is a type ${\rm II}\sb \infty$ factor and $\tau$ is a faithful normal tracial weight on $\M$;
\item $\M=L^\infty[0,\infty)$ and $\tau=\int_0^\infty\, dx$.
\end{itemize}
Note that in each case, $\Aut(\M,\tau)$ acts on $\M$ ergodically.
Recall that $\Aut(\M,\tau)$ acts on $\M$ ergodically if
$\theta(T)=T$ for all $\theta\in \Aut(\M,\tau)$ implies $T=\lambda
1$.  Let $E,F$ be finite projections in $\M$ such that
$\tau(E)=\tau(F)$.  If $\Aut(\M,\tau)$ acts on $\M$ ergodically,
there is a $\theta\in \Aut(\M,\tau)$ such that $\theta(E)=F$.
Furthermore, if $|\!|\!|\cdot|\!|\!|$ is a symmetric gauge semi-norm
on $\J(\M)$, then $|\!|\!|E|\!|\!|=|\!|\!|F|\!|\!|$. In this case, a
semi-norm  $|\!|\!|\cdot|\!|\!|$ on $\J(\M)$ is called a
\emph{normalized} symmetric gauge semi-norm if $|\!|\!|E|\!|\!|=1$
whenever $\tau(E)=1$.\\

Let $(\M,\tau)$ be one of the above semi-finite von Neumann
algebras.
 For every (non-zero) finite projection $E$ in $\M$,
let $\M_E=E\M E$ and $\tau_E(ETE)=\frac{\tau(ETE)}{\tau(E)}$. Then
$(\M_E,\tau_E)$ is a finite von Neumann algebra satisfying the
\emph{weak Dixmier property} (see~\cite{FHNS}), i.e., for every
positive operator $T\in \M_E$, $\tau_E(T)E$ is in the operator
norm closure of the convex hull of $\{S\in \M_E:\,\, \text{$S$ and
$T$ are equi-measurable}\}$. So in the following sections we will
always assume that $(\M,\tau)$ satisfies the following conditions:
\begin{itemize}
\item[\bf A.] $(\M,\tau)$ is a semi-finite von Neumann algebra
such that $\Aut(\M,\tau)$ acts on $\M$ ergodically; \item[\bf B.]
for every non-zero finite projection $E$ in $\M$, $(\M_E,\tau_E)$
is a finite von Neumann algebra satisfying the weak Dixmier
property.
\end{itemize}

With the above assumptions, it is  easy to show that if $E$ is a
finite projection of $\M$, then $\Aut(\M_E,\tau_E)$ acts on $\M_E$
ergodically.\\

A \emph{simple} operator in a semi-finite von Neumann algebra
$(\M,\tau)$ is an operator $T=a_1E_1+\cdots+a_nE_n$, where
$E_1,\cdots,E_n$ are mutually orthogonal projections. The
following lemma is Corollary 3.7 of~\cite{FHNS}.

\begin{Lemma}
 \label{L:simple operators are dense 1}
 Let $\N$ be a finite  von Neumann algebra with a faithful normal
tracial state $\tau_N$. Let $|\!|\!|\cdot|\!|\!|_1$ and
$|\!|\!|\cdot|\!|\!|_2$ be two gauge invariant semi-norms on $\N$.
Then $|\!|\!|\cdot|\!|\!|_1=|\!|\!|\cdot|\!|\!|_2$ on $\N$ if
$|\!|\!|T|\!|\!|_1=|\!|\!|T|\!|\!|_2$ for all positive simple
operators $T\in \N$.
\end{Lemma}

\begin{Lemma}\label{L:local symmetric} Let $(\M,\tau)$ be a
semi-finite von Neumann algebra satisfying the conditions {\bf A}
and {\bf B} and let $|\!|\!|\cdot|\!|\!|$ be  a  symmetric gauge
semi-norm on $\J(\M)$. If $E\in\M$ is a finite projection, then
the restriction of $|\!|\!|\cdot|\!|\!|$ to $(\M_E,\tau_E)$ is
also a symmetric gauge semi-norm on $(\M_E,\tau_E)$.
\end{Lemma}
\begin{proof} It is obvious that the restriction of
$|\!|\!|\cdot|\!|\!|$ to $(\M_E,\tau_E)$ is also a  gauge semi-norm
on $(\M_E,\tau_E)$. Let $\theta\in \Aut(\M_E,\tau_E)$. Define
$|\!|\!|S|\!|\!|_2=|\!|\!|\theta(S)|\!|\!|$ for $S\in \M_E$. We need
to prove $|\!|\!|\cdot|\!|\!|=|\!|\!|\cdot|\!|\!|_2$ on $\M_E$.
  Let $T=a_1E_1+\cdots+a_nE_n$ be a simple positive operator in $\M_E$, where
  $E_1+\cdots+ E_n=E$. Then
$\theta(T)=a_1\theta(E_1)+\cdots+a_n\theta(E_n)$. Since $\theta\in
\Aut(\M_E,\tau_E)$, $\tau(E_k)=\tau(\theta(E_k))$ for $1\leq k\leq
n$. By the assumption of the lemma, $\Aut(\M,\tau)$ acts on $\M$
ergodically. Therefore, there is a $\theta'\in \Aut(\M,\tau)$ such
that $\theta'(E_k)=\theta(E_k)$ for $1\leq k\leq n$. Hence,
$\theta'(T)=\theta(T)$. Since $|\!|\!|\cdot|\!|\!|$ is  a symmetric
gauge semi-norm on $\J(\M)$,
$|\!|\!|T|\!|\!|=|\!|\!|\theta'(T)|\!|\!|=|\!|\!|\theta(T)|\!|\!|=|\!|\!|T|\!|\!|_2$.
By Lemma~\ref{L:simple operators are dense 1},
$|\!|\!|\cdot|\!|\!|=|\!|\!|\cdot|\!|\!|_2$ on $(\M_E,\tau_E)$. This
implies that the restriction of $|\!|\!|\cdot|\!|\!|$ to
$(\M_E,\tau_E)$ is also a symmetric gauge semi-norm on
$(\M_E,\tau_E)$.
\end{proof}

The following lemma is Theorem 3.27 of~\cite{FHNS}.

\begin{Lemma}\label{L:weak dixmier property}
Let $\N$ be a finite  von Neumann algebra with a faithful normal
tracial state $\tau_N$. Then $\N$ satisfies the weak Dixmier
property if and only if $\N$ satisfies one of the following
conditions:
\begin{enumerate}
\item $\N$ is finite dimensional $($hence atomic$)$ and for every
two non-zero minimal projections $E,F\in \N$, $\tau(E)=\tau(F)$,
or equivalently, $(\N,\tau_\N)$ can be identified as a von Neumann
subalgebra of $(M_n(\cc),\tau_n)$ that contains all diagonal
matrices; \item $\N$ is diffuse.
\end{enumerate}
\end{Lemma}

\begin{Corollary}\label{C:weak dixmier property}
Let $(\M,\tau)$ be a semi-finite von Neumann algebra satisfying
the conditions {\bf A} and {\bf B} and let $|\!|\!|\cdot|\!|\!|$
be a normalized symmetric gauge norm on $\M$. If $F$ is a finite
projection in $\M$ such that $\tau(F)\geq 1$, then
$|\!|\!|F|\!|\!|\geq 1$.
\end{Corollary}
\begin{proof} Let $E_1\in \M$ be a finite projection such that
$\tau(E_1)=1$ and $|\!|\!|E_1|\!|\!|=1$. There exists a finite
projection $E\in \M$ such that $E_1, F\leq E$. By
Lemma~\ref{L:local symmetric}, the restriction of
$|\!|\!|\cdot|\!|\!|$ to $(\M_E,\tau_E)$ is also a symmetric gauge
semi-norm on $(\M_E,\tau_E)$. Since $\M_E$ satisfies the weak
Dixmier property, there is a projection $F_1\in \M_E$ such that
$F_1\leq F$ and $\tau(F_1)=1$ by Lemma~\ref{L:weak dixmier
property}. Since $\Aut(\M_E,\tau_E)$ acts on $\M_E$ ergodically,
$|\!|\!|F_1|\!|\!|=|\!|\!|E_1|\!|\!|=1$. By Corollary~\ref{C:S<T},
$|\!|\!|F|\!|\!|\geq |\!|\!|F_1|\!|\!|=1$.
\end{proof}

The following lemma is Corollary 3.36 of~\cite{FHNS}
\begin{Lemma}
\label{T:measure-preserving} Let $\N$ be a finite  von Neumann
algebra with a faithful normal tracial state $\tau_N$. Suppose
$\N$ satisfies the weak Dixmier property and $\Aut(\N,\tau)$ acts
on $\N$ ergodically. If $|\!|\!|\cdot|\!|\!|$ is  a  symmetric
gauge semi-norm on $\N$, then for every $T\in \N$,
\[\|T\|_1\cdot |\!|\!|1|\!|\!|\leq |\!|\!|T|\!|\!|\leq \|T\|\cdot
|\!|\!|1|\!|\!|,
\] where $\|T\|_1=\tau_\N(|T|)$.
\end{Lemma}

\begin{Corollary}\label{C:semi-norm is norm} Let $(\M,\tau)$ be a semi-finite von Neumann algebra satisfying
the  conditions {\bf A} and {\bf B}. If $|\!|\!|\cdot|\!|\!|$ is a
normalized symmetric gauge semi-norm on $\J(\M)$, then
$|\!|\!|\cdot|\!|\!|$ is a symmetric gauge norm on $\J(\M)$.
\end{Corollary}
\begin{proof} Let $T\in \J(\M)$. Then there is a finite projection $E$
in $\M$ such that $T=ETE\in (\M_E,\tau_E)$. We may assume that
$\tau(E)\geq 1$. By Lemma~\ref{L:local symmetric}, the restriction
of $|\!|\!|\cdot|\!|\!|$ to $(\M_E,\tau_E)$ is also a symmetric
gauge semi-norm.  If $T\neq 0$, then by
Lemma~\ref{T:measure-preserving} and Corollary~\ref{C:weak dixmier
property}, $|\!|\!|T|\!|\!|\geq \tau_E(|T|)\cdot
|\!|\!|E|\!|\!|>0$. So $|\!|\!|\cdot|\!|\!|$ is a symmetric gauge
norm on $\J(\M)$.
\end{proof}

The following lemma is Corollary 4.4 of~\cite{FHNS}.
\begin{Lemma}\label{L: type I_n operators are dense}
Let $\N$ be a finite  von Neumann algebra with a faithful normal
tracial state $\tau_N$,  and let $|\!|\!|\cdot|\!|\!|_1$,
$|\!|\!|\cdot|\!|\!|_2$ be two symmetric gauge norms on $\N$.
Suppose $\N$ satisfies the weak Dixmier property and
$\Aut(\N,\tau)$ acts on $\N$ ergodically. Then
$|\!|\!|\cdot|\!|\!|_1=|\!|\!|\cdot|\!|\!|_2$ on $\N$ if
$|\!|\!|T|\!|\!|_1=|\!|\!|T|\!|\!|_2$ for every operator
$T=a_1E_1+\cdots+a_nE_n$ in $\N$ such that $a_1,\cdots,a_n\geq 0$
and $\tau_\N(E_1)=\cdots=\tau_\N(E_n)=\frac{1}{n}$,
$n=1,2,\cdots$.
\end{Lemma}

\begin{Lemma}\label{L:simple operators are
dense}Let $(\M,\tau)$ be a semi-finite von Neumann algebra
satisfying the conditions {\bf A} and {\bf B}. Suppose
$|\!|\!|\cdot|\!|\!|_1$ and $|\!|\!|\cdot|\!|\!|_2$ are two
symmetric gauge norms  on $\J(\M)$. Then
$|\!|\!|\cdot|\!|\!|_1=|\!|\!|\cdot|\!|\!|_2$ on $\J(\M)$ if
$|\!|\!|T|\!|\!|_1=|\!|\!|T|\!|\!|_2$ for every simple positive
operator $T$ in $\J(\M)$ such that $T=a_1E_1+\cdots+a_nE_n$ and
$\tau(E_1)=\cdots=\tau(E_n)$.
\end{Lemma}
\begin{proof} Suppose $|\!|\!|T|\!|\!|_1=|\!|\!|T|\!|\!|_2$ for every
simple operator $T$ in $\J(\M)$. Let $S\in \J(\M)$.  Then there is a
 finite projection $E$ in $\M$ such that $S=ESE\in \M_E$. By
 Lemma~\ref{L:local symmetric}, the restrictions of
 $|\!|\!|\cdot|\!|\!|_1$ and $|\!|\!|\cdot|\!|\!|_2$ to
 $(\M_E,\tau_E)$ are two symmetric gauge norms.
Since $|\!|\!|T|\!|\!|_1=|\!|\!|T|\!|\!|_2$ for every simple
operator $T$ in $\M_E$ such that $T=a_1E_1+\cdots+a_nE_n$ and
$\tau(E_1)=\cdots=\tau(E_n)$,
$|\!|\!|\cdot|\!|\!|_1=|\!|\!|\cdot|\!|\!|_2$ by Lemma~\ref{L: type
I_n operators are dense}.
\end{proof}

\begin{Proposition}\label{P:unitarily invariant implies trace-preserving} Let $(\M,\tau)$ be a semi-finite factor and
let $|\!|\!|\cdot|\!|\!|$ be a norm on $\J(\M)$. Then the
following conditions are equivalent:
\begin{enumerate}
\item $|\!|\!|\cdot|\!|\!|$ is a symmetric gauge norm;
\item $|\!|\!|\cdot|\!|\!|$ is a unitarily invariant norm.
\end{enumerate}
\end{Proposition}
\begin{proof} $``1\Rightarrow 2"$ is obvious. We only prove
$``2\Rightarrow 1"$. We need to prove that for every positive
operator $T\in \J(\M)$ and $\theta\in \Aut(\M,\tau)$,
$|\!|\!|\theta(T)|\!|\!|=|\!|\!|T|\!|\!|$. Let $S=\theta(T)$. Then
$S\in \J(\M)$. Therefor, there is a finite projection $E$ in $\M$
such that $S,T\in \M_E$. By the spectral decomposition theorem,
there is a sequence of simple positive operators $T_n\in \M_E$
such that $S_n=\theta(T_n)\in \M_E$ and
$\lim_{n\rightarrow\infty}\|T_n-T\|=\lim_{n\rightarrow\infty}\|S_n-S\|=0$.
By Lemma~\ref{T:measure-preserving}, $|\!|\!|T-T_n|\!|\!|\leq
\|T-T_n\|\cdot |\!|\!|E|\!|\!|$ and $|\!|\!|S-S_n|\!|\!|\leq
\|S-S_n\|\cdot |\!|\!|E|\!|\!|$. Hence,
$\lim_{n\rightarrow\infty}|\!|\!|T-T_n|\!|\!|=\lim_{n\rightarrow\infty}|\!|\!|S-S_n|\!|\!|=0$.
We need only  prove $|\!|\!|T_n|\!|\!|=|\!|\!|S_n|\!|\!|$ for all
$n=1,2,\cdots$. Suppose $T_n=a_1E_1+\cdots+a_mE_m$.  Then
$S_n=\theta(T_n)=a_1F_1+\cdots+a_mF_m$, where $\theta(E_k)=F_k$
for $1\leq k\leq m$. Since $\theta\in \Aut(\M,\tau)$,
$\tau(E_k)=\tau(F_k)$ for $1\leq k\leq m$.  Since $\M$ is a
factor, there is a unitary operator $U\in\M$ such that
$E_k=UF_kU^*$ for $1\leq k\leq m$. Therefore, $S_n=UT_nU^*$ and
$|\!|\!|T_n|\!|\!|=|\!|\!|S_n|\!|\!|$.
\end{proof}

\subsection{Semi-norms associated to von Neumann algebras}

\begin{Definition}\label{D:seminorm}\emph{ Let $\M$ be a von
Neumann algebra $($not necessarily semi-finite$)$. A
\emph{$($generalized$)$ semi-norm associated to $\M$} is a map
$|\!|\!|\cdot|\!|\!|$ from $\M$ to $[0,\infty]$ satisfying the
following properties:
\begin{enumerate}
\item $|\!|\!|\lambda T|\!|\!|=|\lambda|\cdot |\!|\!|T|\!|\!|$,
\item $|\!|\!|S+T|\!|\!|\leq |\!|\!|S|\!|\!|+|\!|\!|T|\!|\!|$
\end{enumerate}
for all $S,T\in \M$ and $\lambda \in \cc$. To make the definition
nontrivial, we always make the following assumption:
$0<|\!|\!|T|\!|\!|<\infty$ for some non-zero element $T\in \M$.}
\end{Definition}

Let $\I=\{T\in \M:\, |\!|\!|T|\!|\!|<\infty\}$. Then $\I$ is called
\emph{the domain} of the semi-norm $|\!|\!|\cdot|\!|\!|$.

\begin{Definition}\label{D:Symmetric gauge}
\emph{ Let $(\M,\tau)$ be a semi-finite von Neumann algebra. A
semi-norm $|\!|\!|\cdot|\!|\!|$ associated to $\M$ is called
\emph{gauge invariant} if for all $T\in\M$,
$|\!|\!|T|\!|\!|=|\!|\!|\,\, |T|\,\,|\!|\!|$; a semi-norm
$|\!|\!|\cdot|\!|\!|$ associated to $\M$ is \emph{unitarily
invariant} if $|\!|\!|UTV|\!|\!|=|\!|\!|T|\!|\!|$ for all $T\in \M$
and unitary operators $U,V\in \M$;  a semi-norm
 $|\!|\!|\cdot|\!|\!|$ associated to a semi-finite von Neumann
 algebra
$(\M,\tau)$ is called \emph{symmetric} if
\[|\!|\!|\theta(T)|\!|\!|=|\!|\!|T|\!|\!|,\,\,\,\, \forall\, T\in \M,\, \theta\in
\Aut(\M,\tau);
\] a semi-norm $|\!|\!|\cdot|\!|\!|$ associated to $(\M,\tau)$ is called a \emph{symmetric gauge semi-norm} if it
is both symmetric and gauge invariant.}
\end{Definition}

Similar to the proof of Proposition~\ref{P:unitarily invariant norms
on J(M)}, we can prove the following proposition.
\begin{Proposition}\label{P:unitarily invariant norms} Let $|\!|\!|\cdot|\!|\!|$ be a semi-norm associated to $\M$. Then the
following statements are equivalent:
\begin{enumerate}
\item $|\!|\!|\cdot|\!|\!|$ is unitarily invariant;
\item $|\!|\!|\cdot|\!|\!|$ is gauge invariant and unitarily conjugate
invariant, i.e.,  $|\!|\!|UTU^*|\!|\!|=|\!|\!|T|\!|\!|$ for all
$T\in\M$ and unitary operators $U\in \M$;
\item $|\!|\!|\cdot|\!|\!|$ is
gauge invariant and $|\!|\!|T|\!|\!|=|\!|\!|T^*|\!|\!|$ for all
$T\in \M$;
\item for all operators
$T,A,B\in \M$, $|\!|\!|ATB|\!|\!|\leq \|A\|\cdot
|\!|\!|T|\!|\!|\cdot \|B\|$.
\end{enumerate}
\end{Proposition}

\begin{Corollary}\label{C:S<T 2}Let $|\!|\!|\cdot|\!|\!|$ be a semi-norm associated to
$\M$. If $S,T\in\M$ and $0\leq S\leq T$, then $|\!|\!|S|\!|\!|\leq
|\!|\!|T|\!|\!|$.
\end{Corollary}

\begin{Corollary}\label{C:equivalent projections} Let $|\!|\!|\cdot|\!|\!|$
be a unitarily invariant semi-norm associated to $\M$ and $E, F$ be
two equivalent projections in $\M$. Then
$|\!|\!|E|\!|\!|=|\!|\!|F|\!|\!|$.
\end{Corollary}

\begin{Lemma}\label{L:E neq 0} Let $|\!|\!|\cdot|\!|\!|$
be a unitarily invariant semi-norm associated to $\M$ and let
$T\in \M$ be a nonzero element such that $|\!|\!|T|\!|\!|<\infty$.
Then there is a nonzero projection $E$ in $\M$ such that
$|\!|\!|E|\!|\!|<\infty$.
\end{Lemma}
\begin{proof} Since $|\!|\!|\cdot|\!|\!|$ is unitarily invariant, we
may assume $T>0$. By the spectral decomposition theorem, there
exist a $\lambda>0$ and a nonzero projection $E$ in $\M$ such that
$T\geq \lambda E$. By Corollary~\ref{C:S<T 2},
$|\!|\!|E|\!|\!|<\infty$.
\end{proof}

The following theorem shows that, up to a scale $a>0$, the operator
norm $\|\cdot\|$ is the unique unitarily invariant semi-norm
associated to  a type ${\rm III}$ factor.
\begin{Theorem}\label{T:unitarily invariant norms on type III}
 Let $\M$ be a type ${\rm III}$ factor and let $|\!|\!|\cdot|\!|\!|$
be a unitarily invariant semi-norm associated to $\M$. Then there
exists  $a>0$ such that $|\!|\!|\cdot|\!|\!|=a \|\cdot\|$, i.e.,
  $|\!|\!|T|\!|\!|=a \|T\|$ for all $T\in \M$.
\end{Theorem}
\begin{proof} By Lemma~\ref{L:E neq 0}, there is a nonzero projection
$E$ in $\M$ such that $|\!|\!|E|\!|\!|<\infty$. If
$|\!|\!|E|\!|\!|=0$, then $|\!|\!|1|\!|\!|=0$ by
Corollary~\ref{C:equivalent projections}. By
Proposition~\ref{P:unitarily invariant norms}, for every $T$ in
$\M$, $|\!|\!|T|\!|\!|\leq \|T\|\cdot |\!|\!|1|\!|\!|=0$. In our
definition of semi-norm, we assume that $|\!|\!|T|\!|\!|>0$ for
some $T\in \M$. Hence $|\!|\!|E|\!|\!|\neq 0$ for some projection
$E$ in $\M$.  We may assume that $|\!|\!|E|\!|\!|=1$. By
Corollary~\ref{C:equivalent projections}, $|\!|\!|F|\!|\!|=1$ for
every non-zero projection in $\M$. In particular,
$|\!|\!|1|\!|\!|=1$. By Proposition~\ref{P:unitarily invariant
norms}, for every $T$ in $\M$, $|\!|\!|T|\!|\!|\leq \|T\|\cdot
|\!|\!|1|\!|\!|=\|T\|$. On the other hand, let $T\in\M$ be a
positive operator and $\epsilon>0$. By the spectral decomposition
theorem, there is a nonzero projection $F$ in $\M$ such that
$T\geq (\|T\|-\epsilon)F$. By Corollary~\ref{C:S<T 2},
$|\!|\!|T|\!|\!|\geq (\|T\|-\epsilon)\cdot
|\!|\!|F|\!|\!|=\|T\|-\epsilon$. Hence $|\!|\!|T|\!|\!|\geq
\|T\|-\epsilon$. This proves that $|\!|\!|T|\!|\!|=\|T\|$ for
every positive operator $T$ in $\M$ and therefore for every
operator $T$ in $\M$.
\end{proof}

We end this section with the following lemma.
\begin{Lemma}Let $(\M,\tau)$ be a semi-finite von Neumann algebra such that
$\Aut(\M,\tau)$ acts on $\M$ ergodically. If $|\!|\!|\cdot|\!|\!|$
is a normalized symmetric gauge semi-norm associated to $\M$ with
domain $\I$, then $\I\supseteq\J(\M)$ and $|\!|\!|\cdot|\!|\!|$ is a
normalized symmetric gauge semi-norm on $\J(\M)$.
\end{Lemma}
\begin{proof} Let $E$ be a finite projection in $\M$ such that
$\tau(E)=1$. Then $|\!|\!|E|\!|\!|=1$. Suppose that $F$ is a
finite projection in $\M$ such that $n\leq \tau(F)<n+1$. Note that
$\Aut(\M,\tau)$ acts on $\M$ ergodically. By induction, there are
mutually orthogonal finite projections $E_1,E_2,\cdots, E_{n+1}$
in $\M$, $\tau(E_1)=\cdots=\tau(E_{n+1})=1$, such that
$E_1+\cdots+E_n\leq F\leq E_1+\cdots+E_{n+1}$. By
Corollary~\ref{C:S<T 2}, $|\!|\!|F|\!|\!|\leq
|\!|\!|E_1+\cdots+E_{n+1}|\!|\!|\leq n+1$. So every finite
projection is in $\I$. Hence $\I\supseteq \J(\M)$.
\end{proof}

\section{Ky Fan norms associated to semi-finite von Neumann algebras}
Let $(\M,\tau)$ be a semi-finite von Neumann subalgebra of a type
${\rm II}\sb \infty$ factor $(\M_1,\tau_1)$ and let $0\leq t\leq
\infty$. For $T\in \M$, define $|\!|\!|T|\!|\!|_{(t)}$, \emph{the
Ky Fan $t$-th norm of $T$}, by
\[|\!|\!|T|\!|\!|_{(t)}=\left\{\begin{array}{lll}
                                    \|T\|, & \hbox{$t=0$;} \\
                                     \frac{1}{t}\int_0^t \mu_s(T)ds, & \hbox{$0<t\leq
                                     1$;}\\
                                     \int_0^t \mu_s(T)ds, &
                                     \hbox{$1<t\leq\infty$}.
                                   \end{array}
                                 \right.
\]

Let $\U(\M)$ be the set of unitary operators in $\M$ and $\P(\M)$
be the set of projections in $\M$.
\begin{Lemma}\label{L:characterization of Ky fan norms}
 For $0<t\leq 1$, $t|\!|\!|T|\!|\!|_{(t)}=\sup\{|\tau_1(UTE)|:\,
U\in\U(\M_1),\, E\in\P(\M_1),\,  \tau_1(E)=t\}$.
\end{Lemma}
\begin{proof}
We may assume that $T$ is a positive operator. Let $\A$ be a
separable diffuse abelian von Neumann subalgebra of $\M_1$
containing $T$ and let $\alpha$ be a $\ast$-isomorphism from
$(\A,\tau_1)$ onto $(L^\infty[0,\infty), \int_0^\infty\,dx)$ such
that $\tau_1=\int_0^\infty dx\cdot \alpha$. Let $f(x)=\alpha(T)$
and let $f^*(x)$ be the non-increasing rearrangement of $f(x)$.
Then $\mu_s(T)=f^*(s)$. By the definition of $f^*$(see equation
(1) in section 2.1),
\[m(\{f^*>f^*(t)\})=
\lim_{n\rightarrow\infty}m\left(\left\{f^*>f^*(t)+\frac{1}{n}\right\}\right)\leq
t\] and
\[m(\{f^*\geq f^*(t)\})\geq \lim_{n\rightarrow\infty} m\left(\left\{f^*>
f^*(t)-\frac{1}{n}\right\}\right)\geq t.\] Since $f^*$ and $f$ are
equi-measurable, $m(\{f>f^*(t)\})\leq t$ and $m(\{f\geq
f^*(t)\})\geq t$. Therefore, there is a measurable subset $A$ of
$[0,\infty)$, $\{f>f^*(t)\}\subset A\subset \{f\geq f^*(t)\}$, such
that $m(A)=t$. Since $f(x)$ and $f^*(x)$ are equi-measurable,
$\int_A f(s)ds=\int_0^t f^*(s)ds$.  Let $E'=\alpha^{-1}(\chi_A)$.
Then $\tau_1(E')=t$ and $\tau_1(TE')=\int_A f(s)ds=\int_0^t
f^*(s)ds=t|\!|\!|T|\!|\!|_{(t)}$. Hence, $t|\!|\!|T|\!|\!|_{(t)}\leq
\sup\{|\tau_1(UTE)|:\, U\in\U(\M_1), E\in\P(\M_1),\,  \tau_1(E)=t\}$. \\

We need to prove that if $E$ is a projection in $\M_1$,
$\tau_1(E)=t$, and $U\in \U(\M_1)$, then $t|\!|\!|T|\!|\!|_{(t)}\geq
|\tau_1(UTE)|$.
 By the Schwartz
inequality, $|\tau_1(UTE)|=|\tau_1(EUT^{1/2}T^{1/2}E)|\leq
\tau_1(U^*EUT)^{1/2}\tau_1(ET)^{1/2}$. By
Corollary~\ref{C:s-number}, $\tau_1(ET)=\int_0^1 \mu_s(ET)ds$. By
Corollary~\ref{C:property of s-number}, $\mu_s(ET)\leq
\min\{\mu_s(T), \mu_s(E)\|T\|\}$. Note that $\mu_s(E)=0$ for
$s\geq \tau_1(E)=t$. Hence, $\tau_1(ET)\leq
\int_0^t\mu_s(T)ds=t|\!|\!|T|\!|\!|_t$. Similarly,
$\tau_1(U^*EUT)\leq t|\!|\!|T|\!|\!|_t$.  So $|\tau_1(UTE)|\leq
t|\!|\!|T|\!|\!|_t$. This proves that $t|\!|\!|T|\!|\!|_{(t)}\geq
\sup\{|\tau_1(UTE)|:\, U\in\U(\M_1), E\in\P(\M_1),\,
\tau_1(E)=t\}$.

\end{proof}

Similarly, we can prove the following lemma.
\begin{Lemma}\label{L:characterization of Ky fan norms 2}
 For $1\leq t\leq \infty$, $|\!|\!|T|\!|\!|_{(t)}=\sup\{|\tau_1(UTE)|:\,
U\in\U(\M_1),\, E\in\P(\M_1),\,  \tau_1(E)=t\}$.
\end{Lemma}

\begin{Theorem}\label{T:Ky fan norms} For $0\leq t\leq \infty$, $|\!|\!|\cdot|\!|\!|_{(t)}$ is
a normalized
 symmetric gauge norm associated to $(\M,\tau)$.
\end{Theorem}
\begin{proof} By the definition of $s$-number, $\mu_s(T)=
\mu_s(\theta(T))$ for $T\in \M$ and $\theta\in \Aut(\M,\tau)$.
 To prove $|\!|\!|\cdot|\!|\!|_{(t)}$ is a
normalized symmetric gauge norm associated to $(\M,\tau)$, we need
only  prove  the triangle inequality since other parts are
obvious. Let $S, T\in \M$. If $0<t\leq 1$, by
Lemma~\ref{L:characterization of Ky fan norms},
$t|\!|\!|S+T|\!|\!|_{(t)}=\sup\{|\tau_1(U(S+T)E)|:\,
U\in\U(\M_1),\, E\in \P(\M_1),\, \tau_1(E)=t\}\leq
\sup\{|\tau_1(USE)|:\, U\in\U(\M_1),\, E\in \P(\M_1),
\tau_1(E)=t\}+ \sup\{|\tau_1(UTE)|:\, U\in\U(\M_1),\, E\in
\P(\M_1),\,
\tau_1(E)=t\}=t|\!|\!|S|\!|\!|_{(t)}+t|\!|\!|T|\!|\!|_{(t)}$. The
proof of the case $t>1$ is similar.
\end{proof}

The following corollary plays a key role in section 10 (see the
proof of Lemma~\ref{L:embedding into measure topology}).
\begin{Corollary}\label{C:Markov inequality}Let $T\in \M$ and
$\delta>0$. If $|\!|\!|T|\!|\!|_{(1)}<\delta$, then
$\tau(\chi_{(\delta,\infty)}(|T|))\leq
\frac{|\!|\!|T|\!|\!|_{(1)}}{\delta}.$
\end{Corollary}
\begin{proof} We may assume that $\M$ is a type ${\rm II}_\infty$ factor and $T\geq 0$.
 By the proof of Lemma~\ref{L:characterization of Ky fan norms},
\[|\!|\!|T|\!|\!|_{(1)}=\sup\{|\tau(UTE)|:\, U\in\U(\M),\,
E\in\P(\M),\,  \tau(E)\leq 1\}.\] If
$\tau(\chi_{(\delta,\infty)}(T))>1$, then there is a
sub-projection $E$ of $\chi_{(\delta,\infty)}(T)$ such that
$\tau(E)=1$. Then $TE\geq \delta E$. Hence,
$|\!|\!|T|\!|\!|_{(1)}\geq \tau(TE)\geq \tau(\delta E)=\delta$.
This contradicts  the assumption that
$|\!|\!|T|\!|\!|_{(1)}<\delta$. Therefore,
$\tau(\chi_{(\delta,\infty)}(T))\leq 1$. So
 $|\!|\!|T|\!|\!|_{(1)}\geq
\tau(T\chi_{(\delta,\infty)}(T))\geq \tau(\delta
\chi_{(\delta,\infty)}(T))\geq
\delta\tau(\chi_{(\delta,\infty)}(T))$. This implies the corollary.
\end{proof}

\begin{Proposition}\label{P: Monotone property of Ky Fan norms} Let
$(\M,\tau)$ be a semi-finite von Neumann algebra and $T\in
(\M,\tau)$. Then $|\!|\!|T|\!|\!|_{(t)}$ is a non-increasing
continuous function on $[0,1]$ and a non-decreasing continuous
function  on $[1,\infty]$.
\end{Proposition}
\begin{proof} Let $0< t_1<t_2\leq 1$.
$|\!|\!|T|\!|\!|_{(t_1)}-|\!|\!|T|\!|\!|_{(t_2)}=\frac{1}{t_1}\int_0^{t_1}
\mu_s(T)ds- \frac{1}{t_2}\int_0^{t_2}
\mu_s(T)ds=\frac{\frac{1}{t_1}\int_0^{t_1}
\mu_s(T)ds-\frac{1}{t_2-t_1}\int_{t_1}^{t_2}\mu_s(T)ds}{t_2(t_2-t_1)}\leq
0.$ Since $\mu_s(T)$ is right-continuous, $|\!|\!|T|\!|\!|_{(t)}$ is
a non-increasing continuous function on $[0,1]$. Since $\mu_s(T)\geq
0$ for $s\in [0,\infty)$,  $|\!|\!|T|\!|\!|_{(t)}$ is a
non-decreasing continuous function on $[1,\infty]$.
\end{proof}

\begin{Proposition}\label{P:comparison proposition} Let $(\M,\tau)$
be a semi-finite von Neumann algebra satisfying the conditions
{\bf A} and {\bf B} in section 3.4,
  and let
$|\!|\!|\cdot|\!|\!|$ be a normalized symmetric gauge norm on
$\J(\M)$. Then for every $T\in \J(\M)$,
\[|\!|\!|T|\!|\!|_{(1)}\leq |\!|\!|T|\!|\!|.
\]
\end{Proposition}
\begin{proof} We can assume that $T$ is a positive operator in
$\J(\M)$. Then there is a  finite projection $F$ in $\M$ such that
$T=FTF\in \M_F$. We can assume that $\tau(F)=k$ is a positive
integer. By the assumption of the proposition, $(\M_F,\tau_F)$
satisfies the weak Dixmier property. By Lemma~\ref{L:weak dixmier
property}, either $(\M_F,\tau_F)$ is a diffuse von Neumann algebra
or $(\M_F,\tau_F)$ is $\ast$-isomorphic to a von Neumann
subalgebra of $(M_n(\cc), \tau_n)$ that contains all diagonal
matrices. In either case, there is a projection $E$ in $\M$,
$E\leq F$, such that $\tau(E)=1$ and
$|\!|\!|T|\!|\!|_{(1)}=|\!|\!|ETE|\!|\!|_{(1)}$. By
Corollary~\ref{C:symmetric implies unitarily invariant on J(M)}
and Proposition~\ref{P:unitarily invariant norms on J(M)},
$|\!|\!|ETE|\!|\!|\leq |\!|\!|T|\!|\!|$. By
Lemma~\ref{T:measure-preserving}, $|\!|\!|ETE|\!|\!|_{(1)}\leq
|\!|\!|ETE|\!|\!|\leq |\!|\!|T|\!|\!|$.
\end{proof}

\begin{Example}\label{E:Ky fan norms for MnC}\emph{ The Ky Fan $n$-th  norm
of a compact operator $T\in(\B(\H), {\rm Tr})$ is
\[|\!|\!|T|\!|\!|_{({n})}=s_1(T)+\cdots+s_n(T)
\] and
\[|\!|\!|T|\!|\!|_{({\infty})}=s_1(T)+s_2(T)+\cdots.
\]}
\end{Example}

\begin{Corollary}\label{C:comparison corollary} Let
$|\!|\!|\cdot|\!|\!|$ be a normalized unitarily invariant norm on
$\B(\H)$. Then for every $T\in \J(\H)$,
\[s_1(T)\leq |\!|\!|T|\!|\!|\leq s_1(T)+s_2(T)+\cdots.
\]
\end{Corollary}
\begin{proof} By Proposition~\ref{P:comparison proposition},
$s_1(T)=|\!|\!|T|\!|\!|_{(1)}\leq |\!|\!|T|\!|\!|$. On the other
hand, we may assume that $T$ is a positive operator in $\J(\H)$.
Then $T$ is unitarily equivalent to a diagonal operator
$s_1(T)E_1+\cdots+s_n(T)E_n$. Hence,
$|\!|\!|T|\!|\!|=|\!|\!|s_1(T)E_1+\cdots+s_n(T)E_n|\!|\!|\leq
s_1(T)+\cdots+s_n(T)$.
\end{proof}

\section{Dual norms of symmetric gauge norms
on $\J(\M)$}

Throughout this section, we assume that $(\M,\tau)$ is a semi-finite
von Neumann algebra
 satisfying the conditions {\bf A} and {\bf B} in section 3.4.
Recall that $\J(\M)$ is the subset of $\M$ consisting of operators
$T$ in $\M$ such that $T=ETE$ for some finite projection $E\in
\M$. Note that for  two operators $S,T$ in $\J(\M)$, there is a
finite projection $F$ in $\M$ such that $S,T\in \M_F=F\M F$.

\subsection{Dual norms}
Let $|\!|\!|\cdot|\!|\!|$ be a norm on $\J(\M)$. For $T\in \J(\M)$,
define
\[|\!|\!|T|\!|\!|^\#_{\M,\tau}=\sup\{|\tau(TX)|:\, X\in\J(\M),\, |\!|\!|X|\!|\!|\leq 1\}.
\]
When no confusion arises, we simply write $|\!|\!|\cdot|\!|\!|^\#$
or $|\!|\!|\cdot|\!|\!|_\M^\#$ instead of
$|\!|\!|\cdot|\!|\!|^\#_{\M,\tau}$.

\begin{Lemma} $|\!|\!|\cdot|\!|\!|^\#$ is a norm on
$\J(\M)$.
\end{Lemma}
\begin{proof} Note that if $T\in \J(\M)$ is not 0, then
$|\!|\!|T|\!|\!|^\#\geq \frac{\tau(TT^*)}{|\!|\!|T^*|\!|\!|}>0$.
It is easy to check that $|\!|\!|\cdot|\!|\!|^\#$ satisfies other
conditions for a norm.
\end{proof}

\begin{Definition}\label{D:dual norm}
\emph{ $|\!|\!|\cdot|\!|\!|^\#$ is called the \emph{dual norm} of
$|\!|\!|\cdot|\!|\!|$ on $\J(\M)$ with respect to $\tau$.}
\end{Definition}

The following lemma follows simply from the definition of dual norm.
\begin{Lemma}\label{L:Holder inequality for bounded operators}
Let $|\!|\!|\cdot|\!|\!|$ be a norm on $\J(\M)$ and
$|\!|\!|\cdot|\!|\!|^\#$ be the dual norm on $\J(\M)$. Then for
$S,T\in \J(\M)$, $|\tau(ST)|\leq |\!|\!|S|\!|\!|\cdot
|\!|\!|T|\!|\!|^\#.$
\end{Lemma}

For $T\in \M$, define $\|T\|_1=\tau(|T|)$. Then
$\|T\|_1=|\!|\!|T|\!|\!|_{(\infty)}$.  The following corollary is
the H$\ddot{\text{o}}$lder's inequality for operators in $\J(\M)$.
\begin{Corollary}\label{C:Holder inequality for bounded operators}
 Let  $|\!|\!|\cdot|\!|\!|$ be a gauge invariant norm on
$\J(\M)$ and $|\!|\!|\cdot|\!|\!|^\#$ be the dual norm. Then for
$S,T\in \J(\M)$, $\|ST\|_1\leq |\!|\!|S|\!|\!|\cdot
|\!|\!|T|\!|\!|^\#$.
\end{Corollary}
\begin{proof} By Lemma~\ref{L:characterization of Ky fan norms 2},
 $\|ST\|_1=|\!|\!|ST|\!|\!|_{(\infty)}=\sup\{|\tau(UST)|:\,
U\in \U(\M)\}$. By Lemma~\ref{L:Holder inequality for bounded
operators} and Lemma~\ref{L:gauge=left unitarily invariant},
$|\tau(UST)|\leq |\!|\!|US|\!|\!|\cdot
|\!|\!|T|\!|\!|^\#=|\!|\!|S|\!|\!|\cdot |\!|\!|T|\!|\!|^\#.$
\end{proof}

Let $E$ be a (non-zero) finite projection in $\M$. Recall that
$\M_E=E\M E$ is a finite von Neumann algebra with a faithful normal
tracial state $\tau_E$ such that $\tau_E(T)=\frac{\tau(T)}{\tau(E)}$
for $T\in \M_E$. If  $|\!|\!|\cdot|\!|\!|$ is a norm on $\M_E$, the
dual norm of $T\in \M_E$ with respect to $\tau_E$ is defined by
\[|\!|\!|T|\!|\!|^\#_{\M_E,\,\tau_E}=\sup\{|\tau_E(TX)|:\, X\in\M_E,\, |\!|\!|X|\!|\!|\leq 1\}.
\]
\begin{Lemma}\label{L:reduction to finite case} Suppose $|\!|\!|\cdot|\!|\!|$ is
a unitarily invariant norm on $\J(\M)$. Let $E$ be a non-zero finite
projection in $\M$ and $T\in \M_E$. Then
\[|\!|\!|T|\!|\!|^\#_{\M,\tau}=\tau(E)\cdot|\!|\!|T|\!|\!|^\#_{\M_E,\,\tau_E}. \]
\end{Lemma}
\begin{proof}Since $T=ETE$, for every $X\in \J(\M)$,
$\tau(TX)=\tau(ETEX)=\tau(ETEEXE)=\tau(E)\cdot \tau_E(ETEEXE)$. If
$|\!|\!|X|\!|\!|\leq 1$, then $|\!|\!|EXE|\!|\!|\leq
|\!|\!|X|\!|\!|$ by Proposition~\ref{P:unitarily invariant norms on
J(M)}. This implies that
\[|\!|\!|T|\!|\!|^\#_{\M,\tau}=\sup\{|\tau(TX)|:\, X\in\J(\M),\, |\!|\!|X|\!|\!|\leq
1\}=\sup\{|\tau(TX)|:\, X\in\M_E,\, |\!|\!|X|\!|\!|\leq 1\}\]
\[=\tau(E)\cdot\sup\{|\tau_E(TX)|:\, X\in\M_E,\, |\!|\!|X|\!|\!|\leq
1\}=\tau(E)\cdot|\!|\!|T|\!|\!|^\#_{\M_E,\,\tau_E}.
\]
\end{proof}

The next lemma follows from Proposition 6.5, Proposition 6.6 and
Theorem 6.10 of~\cite{FHNS}.

\begin{Lemma}\label{L:dual invariant and unitarily norms}
Let $\N$ be a finite von Neumann algebra with a faithful normal
tracial state $\tau_\N$. We have the following:
 \begin{enumerate}
\item if $|\!|\!|\cdot|\!|\!|$ is a unitarily invariant norm on
$\N$, then
 $|\!|\!|\cdot|\!|\!|_{\N,\tau_\N}^\#$ is  also a unitarily invariant norm on
$\N$; \item if $|\!|\!|\cdot|\!|\!|$ is a symmetric gauge norm on
$\N$, then
 $|\!|\!|\cdot|\!|\!|_{\N,\tau_\N}^\#$ is also a symmetric gauge norm on
$\N$. Furthermore, if $|\!|\!|1|\!|\!|=1$, then
$|\!|\!|1|\!|\!|_{\N,\tau_\N}^\#=1$.
 \end{enumerate}
\end{Lemma}

 Combining
Lemma~\ref{L:reduction to finite case} and Lemma~\ref{L:dual
invariant and unitarily norms}, we have the following proposition.

\begin{Proposition}\label{P:dual invaraint norm}
 Let  $|\!|\!|\cdot|\!|\!|$ be a  norm on $\J(\M)$. We have the
 following:
 \begin{enumerate}
\item if $|\!|\!|\cdot|\!|\!|$ is a unitarily invariant norm on $\J(\M)$, then
 $|\!|\!|\cdot|\!|\!|^\#$ is  also a unitarily invariant norm on
$\J(\M)$;
\item if $|\!|\!|\cdot|\!|\!|$ is a symmetric gauge norm on $\J(\M)$, then
 $|\!|\!|\cdot|\!|\!|^\#$ is also a symmetric gauge norm on
$\J(\M)$. Furthermore, if $|\!|\!|\cdot|\!|\!|$ is a
 normalized  norm, i.e., $|\!|\!|E|\!|\!|=1$ if
 $\tau(E)=1$, then $|\!|\!|\cdot|\!|\!|^\#$ is also a normalized
  norm.
 \end{enumerate}
\end{Proposition}

The next Lemma follows from Lemma 6.9 of~\cite{FHNS}.

\begin{Lemma}\label{L:Dual norm of simple operators in finite vNA}
 Let $\N$ be a finite von Neumann algebra with a faithful normal
 tracial state $\tau_\N$. Suppose $\N$
satisfies the weak Dixmier property and $|\!|\!|\cdot|\!|\!|$ is a
symmetric gauge norm on $\N$. If $T=a_1 E_1+\cdots+a_n E_n$ is a
positive simple operator in $\N$, then
\begin{equation*}
|\!|\!|T|\!|\!|_{\N,\tau_\N}^\#=\sup\left\{\sum_{k=1}^n
a_kb_k\tau(E_k):\,\, X=b_1E_1+\cdots+b_nE_n \geq 0 \text{ and }
|\!|\!|X|\!|\!|\leq 1\right\}.
\end{equation*}
\end{Lemma}

\begin{Lemma}\label{L:Dual norm of simple operators} Let $|\!|\!|\cdot|\!|\!|$ be a symmetric gauge norm on $\J(\M)$.
 If $T=a_1 E_1+\cdots+a_n E_n$ is a positive simple
operator in $\J(\M)$, then
\[|\!|\!|T|\!|\!|^\#=\sup\left\{\sum_{k=1}^n a_kb_k\tau(E_k):\,\,
S=b_1E_1+\cdots+b_nE_n\geq 0,\,|\!|\!|S|\!|\!|\leq 1\right\}.
\]
\end{Lemma}
\begin{proof} Let $E=E_1+\cdots+E_n$.
 By Lemma~\ref{L:reduction to finite case} and Lemma~\ref{L:Dual norm of simple operators in finite vNA},
\begin{eqnarray*}
   |\!|\!|T|\!|\!|^\#&=& \tau(E)\cdot|\!|\!|T|\!|\!|^\#_{\M_E,\,\tau_E} \\
   &=& \tau(E)\sup\left\{\sum_{k=1}^n a_kb_k\tau_E(E_k):\,\,
S=b_1E_1+\cdots+b_nE_n\geq 0,\,|\!|\!|S|\!|\!|\leq 1\right\} \\
   &=& \sup\left\{\sum_{k=1}^n a_kb_k\tau(E_k):\,\,
S=b_1E_1+\cdots+b_nE_n\geq 0,\,|\!|\!|S|\!|\!|\leq 1\right\}.
\end{eqnarray*}

\end{proof}

\subsection{Dual norms of Ky Fan norms}

The following lemma is Lemma 6.11 of~\cite{FHNS}.
\begin{Lemma}\label{L:dual of ky fan norms for matrix} For $T\in (M_n(\cc),\tau_n)$,
\[|\!|\!|T|\!|\!|_{(\frac{k}{n}),\, M_n(\cc),\, \tau_n}^\#=\max\left\{\frac{k}{n}\|T\|,
\|T\|_{1,\,\tau_n}\right\},
\] where
$|\!|\!|T|\!|\!|_{(\frac{k}{n}),\, M_n(\cc),\,
\tau_n}=\frac{s_1(T)+\cdots+s_k(T)}{k}$ and
$\|T\|_{1,\,\tau_n}=\tau_n(|T|)=\frac{s_1(T)+\cdots+s_n(T)}{n}$.
\end{Lemma}

\begin{Theorem}\label{T:dual of ky fan norms for matrix} For $T\in
\J(\H)$ and $k=1,2,\cdots,\infty$,
\[|\!|\!|T|\!|\!|_{(k)}^\#=\max\left\{\|T\|,
\frac{1}{k}\|T\|_1\right\},
\] where
$|\!|\!|T|\!|\!|_{(k)}=s_1(T)+\cdots+s_k(T)$, $\|T\|_1={\rm
Tr}(|T|)=s_1(T)+s_2(T)+\cdots$ and $\frac{1}{\infty}=0$.
\end{Theorem}
\begin{proof}  For $T\in \J(\H)$, there is a finite rank
projection $E$ such that $T=ETE\in \B(\H)_E$. Let ${\rm Tr}(E)=n$.
Then $\B(\H)_E\cong M_n(\cc)$. First assume $k<\infty$. We may
assume that $n\geq k$. Then
$|\!|\!|T|\!|\!|_{(k)}=k|\!|\!|T|\!|\!|_{(\frac{k}{n}),\,\tau_n}$.
By Lemma~\ref{L:reduction to finite case} and Lemma~\ref{L:dual of
ky fan norms for matrix},
\[|\!|\!|T|\!|\!|_{(k)}^\#={\rm Tr}(E)\cdot\left(k|\!|\!|T|\!|\!|_{(\frac{k}{n}),\, M_n(\cc),\, \tau_n}^\#\right)=
\frac{n}{k}\max\left\{\frac{k}{n}\|T\|,
\|T\|_{1,\,\tau_n}\right\}=\max\left\{\|T\|,
\frac{1}{k}\|T\|_{1}\right\}.\] If $k=\infty$, then
$|\!|\!|T|\!|\!|_{(\infty)}^\#=|\!|\!|T|\!|\!|_{(n)}^\#$ by
Lemma~\ref{L:Dual norm of simple operators}.  Since
$\frac{1}{n}\|T\|_1\leq \|T\|$,
$|\!|\!|T|\!|\!|_{(\infty)}^\#=|\!|\!|T|\!|\!|_{(n)}^\#=\max\left\{\|T\|,
\frac{1}{n}\|T\|_1\right\}=\|T\|$.
\end{proof}

The following Lemma is Theorem 6.19 of~\cite{FHNS}.
\begin{Lemma}\label{T:dual norms of Ky Fan norms} Let $\N$ be a
type ${\rm II}\sb 1$ factor with the faithful normal tracial state
$\tau_\N$ and $0\leq t\leq 1$. Then
\[(|\!|\!|T|\!|\!|_{(t),\,\N,\,\tau_\N}^\#=\max\{t\|T\|,\|T\|_{1,\,\tau_\N}\},\quad\forall\, T\in\M,
\]
where $|\!|\!|T|\!|\!|_{(t),\,\N,\,\tau_\N}$ is the Ky Fan $t$-th
norm of $T$ with respect to $\tau_\N$ and
$\|T\|_{1,\,\tau_\N}=\tau_\N(|T|)$.
\end{Lemma}

\begin{Theorem}\label{T:dual norms of Ky Fan norms for II infty} Let $\M$ be a
type ${\rm II}\sb \infty$ factor and $0\leq t\leq \infty$. Then for
all $T\in \J(\M)$, \[|\!|\!|T|\!|\!|_{(t)}^\#=\left\{
  \begin{array}{ll}
    \max\{t\|T\|,\|T\|_1\}, & \hbox{if $0\leq t\leq 1$;} \\
    \max\{\|T\|,\frac{1}{t}\|T\|_1\}, & \hbox{if $1<t\leq \infty$.}
  \end{array}
\right.\]
\end{Theorem}
\begin{proof} Let $T\in \J(\M)$ and $0<t<\infty$.  There is a finite projection $E$
in $\M$ such that $T=ETE$ is in $\M_{E}$. We can assume that
$\tau(E)=n>t$. Let $\tau_E(ESE)=\frac{\tau(ESE)}{\tau(E)}$. Then
$(\M_E,\tau_E)$ is a type ${\rm II}\sb 1$ factor and $\tau_E$ is the
unique tracial state on $\M_E$. If $0<t\leq 1$, by
Lemma~\ref{L:characterization of Ky fan norms},
\begin{align*}
t|\!|\!|T|\!|\!|_{(t)}&=\sup\{|\tau(UTE')|:\, U\in\U(\M_E),\,
E'\in\P(\M_E),\,  \tau(E')=t\}\\
& =\tau(E)\cdot \sup\{|\tau_E(UTE')|:\, U\in\U(\M_E),\,
E'\in\P(\M_E),\,
\tau_E(E')=t/n\}\\
&=\tau(E)\frac{t}{n}|\!|\!|T|\!|\!|_{(\frac{t}{n}),\,\M_E,\,\tau_E}=t|\!|\!|T|\!|\!|_{(\frac{t}{n}),\,\M_E,\,\tau_E},
\end{align*} where
$|\!|\!|T|\!|\!|_{(\frac{t}{n}),\,\M_E,\,\tau_E}$ means the
 Ky Fan $\frac{t}{n}$-th norm of $T\in \M_E$ with respect to the
tracial state $\tau_E$.\\

Hence,
$|\!|\!|T|\!|\!|_{(t)}=|\!|\!|T|\!|\!|_{(\frac{t}{n}),\,\M_E,\,\tau_E}$.
By Lemma~\ref{L:reduction to finite case} and Lemma~\ref{T:dual
norms of Ky Fan norms for II infty},
\[|\!|\!|T|\!|\!|_{(t)}^\#=\tau(E)\cdot\left(|\!|\!|T|\!|\!|_{(\frac{t}{n}),\,\M_E,\,\tau_E}^\#\right)=
n\max\left\{\frac{t}{n}\|T\|,\|T\|_{1,\,\tau_E}\right\}=\max\{t\|T\|,\|T\|_1\}.\]
If $1<t<\infty$, then
$|\!|\!|T|\!|\!|_{(t)}=t|\!|\!|T|\!|\!|_{(\frac{t}{n}),\,\M_E,\,\tau_E}$.
By Lemma~\ref{L:reduction to finite case} and Lemma~\ref{T:dual
norms of Ky Fan norms},
\[|\!|\!|T|\!|\!|_{(t)}^\#=\tau(E)\cdot\left(t|\!|\!|T|\!|\!|_{(\frac{t}{n}),\,\M_E,\,\tau_E}^\#\right)=
\frac{n}{t}\max\left\{\frac{t}{n}\|T\|,\|T\|_{1,\,\tau_E}\right\}=\max\left\{\|T\|,\frac{1}{t}\|T\|_1\right\}.\]
Similar to the proof of Theorem~\ref{T:dual of ky fan norms for
matrix}, $|\!|\!|T|\!|\!|_{(\infty)}^\#=\|T\|$.
\end{proof}
\subsection{Second dual norms}

The following lemma follows from Theorem~{\bf C} of~\cite{FHNS}.
\begin{Lemma}\label{L:second dual of finite vNA} Let $(\N,\tau_\N)$ be a finite von
Neumann algebra satisfying the weak Dixmier property and let
$|\!|\!|\cdot|\!|\!|$ be a symmetric gauge norm on $\M$. Then
 $|\!|\!|\cdot|\!|\!|^{\#\#}=|\!|\!|\cdot|\!|\!|$.
\end{Lemma}

\begin{Theorem}\label{T:Theorem B}
Let $(\M,\tau)$ be a semi-finite von Neumann algebra satisfying
the conditions {\bf A} and {\bf B} in section 3.4.
 If $|\!|\!|\cdot|\!|\!|$ is a symmetric gauge norm on $\J(\M)$,
then $|\!|\!|\cdot|\!|\!|^\#$ is also a symmetric gauge norm on
$\J(\M)$ and $|\!|\!|\cdot|\!|\!|^{\#\#}=|\!|\!|\cdot|\!|\!|$ on
$\J(\M)$.
\end{Theorem}
\begin{proof}By Proposition~\ref{P:dual
invaraint norm}, $|\!|\!|\cdot|\!|\!|^{\#}$ is a symmetric gauge
norm on $\J(\M)$. Furthermore,
  both $|\!|\!|\cdot|\!|\!|^{\#\#}$ and
$|\!|\!|\cdot|\!|\!|$ are symmetric gauge norms on $\J(\M)$. We
need to prove that $|\!|\!|T|\!|\!|=|\!|\!|T|\!|\!|^{\#\#}$ for
every positive operator $T\in\J(\M)$. Let $E$ be a finite
projection in $\M$ such that $T\in \M_E$.  By
Lemma~\ref{L:reduction to finite case} and Lemma~\ref{L:second
dual of finite vNA},
\begin{eqnarray}
  |\!|\!|T|\!|\!|^{\#\#}_{\M,\tau} &=& \sup\{|\tau(TX)|:\, X\in\J(\M),\, |\!|\!|X|\!|\!|^\#_{\M,\tau}\leq
1\}\nonumber \\
   &=& \sup\{\tau(E)\cdot |\tau_E(TX)|:\, X\in\M_E,\, |\!|\!|X|\!|\!|^\#_{\M,\tau}\leq
1\}\nonumber \\
   &=& \sup\{|\tau_E(T(\tau(E)X))|:\, X\in\M_E,\, |\!|\!|\tau(E)X|\!|\!|^\#_{\M_E,\,\tau_E}\leq
1\}\nonumber \\
   &=&|\!|\!|T|\!|\!|^{\#\#}_{\M_E,\,\tau_E}=|\!|\!|T|\!|\!|_{\M_E,\,\tau_E}=|\!|\!|T|\!|\!|.\nonumber
\end{eqnarray}
\end{proof}

\section{Main result}

 Throughout this section, we assume that $(\M,\tau)$ is a
semi-finite von Neumann algebra.

\begin{Lemma}\label{L:T_f} Let $f(x)=\sum_{k=1}^n a_k
\chi_{[\alpha_{k-1},\,\alpha_k)}(x)$, where $a_1\geq a_2\geq\cdots
\geq a_n\geq 0(=a_{n+1})$ and
$0=\alpha_0<\alpha_1<\cdots<\alpha_n<\infty$.
 For
$T\in \M$, define
\[|\!|\!|T|\!|\!|_{f}=\int_0^\infty f(s)\mu_s(T)ds.
\] Then
\[|\!|\!|T|\!|\!|_{f}=\sum_{k=1}^n \min\{\alpha_k,1\}(a_k-a_{k+1})|\!|\!|T|\!|\!|_{(\alpha_k)}.\]
\end{Lemma}
\begin{proof}  Since
$t|\!|\!|T|\!|\!|_{(t)}=\int_{0}^t \mu_s(T)ds$ for $0\leq t\leq 1$
and $|\!|\!|T|\!|\!|_{(t)}=\int_{0}^t \mu_s(T)ds$ for $1\leq
s<\infty$, summation by parts shows that
\begin{align*}
 |\!|\!|T|\!|\!|_{f}&=\int_0^\infty f(s)\mu_s(T)ds=
  a_1\int_0^{\alpha_1}\mu_s(T)ds+a_2\int_{\alpha_1}^{\alpha_2}\mu_s(T)ds+
   \cdots+a_n\int_{\alpha_{n-1}}^{\alpha_n}\mu_s(T)ds\\
  &=\sum_{k=1}^n \min\{\alpha_k,1\}(a_k-a_{k+1})|\!|\!|T|\!|\!|_{(\alpha_k)}.
  \end{align*}
\end{proof}
\begin{Corollary}\label{L:T_f is a norm} $|\!|\!|\cdot|\!|\!|_{f}$  is a
symmetric  gauge norm associated   to $\M$ and therefore a
symmetric gauge norm on $\J(\M)$. Furthermore, if $\tau(E)=1$ then
$|\!|\!|E|\!|\!|_{f}=\int_0^1 f(x)dx$.
\end{Corollary}

\begin{Lemma}\label{L:proof of theorem c} Let $(\M,\tau)$ be a semi-finite von Neumann algebra
and $E\in \M$ be a $($non-zero$)$ finite projection. Suppose $\M_E$
is a diffuse von Neumann algebra and $T,X\in \M_E$ are positive
operators such that $T=a_1E_1+\cdots+a_nE_n$, $E_1+\cdots+E_n=E$,
and $\tau(E_1)=\cdots=\tau(E_n)$. Then there is  a sequence of
simple positive operators $X_n\in\M_E$ satisfying the following
conditions:
\begin{enumerate}
\item $0\leq X_1\leq X_2\leq \cdots\leq X$ and hence $0\leq
\mu_s(X_1)\leq \mu_s(X_2)\leq \cdots\leq \mu_s(X)$ for all $s\in
[0,\infty)$; \item $\lim_{n\rightarrow \infty}\mu_s(X_n)=\mu_s(X)$
for almost all $s\in [0,\infty)$; \item there exists an $r_n\in
\nn$ such that
$T={a}_{n,1}{E}_{n,1}+\cdots+{a}_{n,r_n}{E}_{n,r_n}$ and
$X_n=b_{n,1}F_{n,1}+\cdots+b_{n,r_n}F_{n,r_n}$, where
$E_{n,1}+\cdots+E_{n,r_n}=F_{n,1}+\cdots+F_{n,r_n}=E$ and
$\tau(E_{n,i})=\tau(F_{n,j})$ for $1\leq i,j\leq r_n$.
\end{enumerate}
\end{Lemma}
\begin{proof} Since $\M_E$ is diffuse, there is a separable diffuse
abelian von Neumann subalgebra $\A$ of $\M_E$ such that $X\in \A$.
Let $\theta$ be a $\ast$-isomorphism  from $\A$ onto
$L^\infty[0,1]$ such that $\tau_E=\int_0^1 dx\cdot\theta$. Let
$f(x)=\theta(X)$. We can choose a sequence of simple functions
$f_n(x)$ in $L^\infty[0,1]$ such that $0\leq f_1(x)\leq f_2(x)\leq
\cdots \leq f(x)$ and $\lim_{n\rightarrow\infty}f_n(x)=f(x)$ for
almost all $x$. Let $X_n=\theta^{-1}(f_n(x))$. Then $X_n\in \M_E$
and $0\leq X_1\leq X_2\leq \cdots\leq X$. By
Lemma~\ref{L:s-number2},
\[\mu_s(T)=\inf\{\|TF\|:\, F\in \P(\M),\quad \tau(F^\perp)=s\}\]
\[=\inf\{\|TF\|:\, F\in \P(\M_E),\quad
\tau_E(F^\perp)=s\tau(E)\}=f^*(\tau(E)s),
\] where $f^*(x)$ is the non-increasing rearrangement of $f(x)$.
  Therefore, we obtain 1 and 2. To obtain 3, we need only
construct
$f_n(x)=\alpha_{n,1}\chi_{I_{n,1}}(x)+\cdots+\alpha_{n,r_n}\chi_{I_{n,r_n}}(x)$
such that $m(I_{n,1})=\cdots=m(I_{n,r_n})=\frac{\tau_E(E_1)}{k_n}$,
for some $k_n\in \nn$.
\end{proof}

Let $\F$ be the set of non-increasing, non-negative, right
continuous simple functions $f(x)$ on $[0,\infty)$  with compact
supports  such that $\int_0^1 f(x)dx\leq 1$. $\forall\, f(x)\in
\F$, $f(x)=\sum_{k=1}^n a_k \chi_{[\alpha_{k-1},\,\alpha_k)}(x)$,
where $a_1\geq a_2\geq\cdots \geq a_n\geq 0(=a_{n+1})$ and
$0=\alpha_0<\alpha_1<\cdots<\alpha_n<\infty$.\\

 Recall that a normalized
norm $|\!|\!|\cdot|\!|\!|$ on $\J(\M)$ of a semi-finite von
Neumann algebra $\M$ is a norm on $\J(\M)$  such that
$|\!|\!|E|\!|\!|=1$ for some projection $E$ with $\tau(E)=1$. The
following theorem is the main result of this paper.\\

\begin{Theorem}\label{T:Theorem C}
Let $(\M,\tau)$ be a
 semi-finite von Neumann algebra satisfying the
 conditions {\bf A} and {\bf B} in section 3.4.
If  $|\!|\!|\cdot|\!|\!|$ is a normalized symmetric gauge norm on
$\J(\M)$, then  there is a subset $\F'$ of $\F$ containing the
characteristic function on $[0,1]$ such that for all $T\in
\J(\M)$, $|\!|\!|T|\!|\!|=\sup\{|\!|\!|T|\!|\!|_{f}:\, f\in
\F'\}$.
\end{Theorem}
\begin{proof} Suppose $|\!|\!|\cdot|\!|\!|$ is
a normalized symmetric gauge norm on $\M$. Let $\F'=\{\mu_s(X):\, X
\,\text{is a simple positive operator in}\, \J(\M),\,
|\!|\!|X|\!|\!|^\#\leq 1\}$.
 For every  positive operator $X\in \J(\M)$ such that $|\!|\!|X|\!|\!|^\#\leq 1$,
  by Proposition~\ref{P:comparison proposition}, $\int_0^1\mu_s(X)ds=|\!|\!|X|\!|\!|_{(1)}\leq
|\!|\!|X|\!|\!|^\#\leq 1$. If $E$ is a  projection such that
$\tau(E)=1$, then $|\!|\!|E|\!|\!|^\#=1$ by
Proposition~\ref{P:dual invaraint norm}. Note that
$\mu_s(E)=\chi_{[0,1]}(s)$.  Therefore, $\F'\subset \F$ and
$\chi_{[0,1]}(s)\in \F'$. For $T\in \J(\M)$, define
\[|\!|\!|T|\!|\!|'=\sup\{|\!|\!|T|\!|\!|_{f}:\, f\in \F'\}.
\] By Corollary~\ref{L:T_f is a norm}, $|\!|\!|\cdot|\!|\!|'$ is a symmetric gauge norm on
$\J(\M)$.
 By Lemma~\ref{L:simple operators are dense}, to prove that
$|\!|\!|\cdot|\!|\!|'=|\!|\!|\cdot|\!|\!|$,  we need to prove
$|\!|\!|T|\!|\!|'=|\!|\!|T|\!|\!|$ for every positive simple
operator $T\in \J(\M)$ such that $T=a_1E_1+\cdots+a_nE_n$ and
$\tau(E_1)=\cdots=\tau(E_n)=c>0$. \\

 By
Lemma~\ref{L:Dual norm of simple operators} and
Theorem~\ref{T:Theorem B},
\[|\!|\!|T|\!|\!|=\sup\left\{c\sum_{k=1}^n a_kb_k:\,\,
X=b_1E_1+\cdots+b_nE_n\geq 0,\,|\!|\!|X|\!|\!|^\#\leq 1\right\}.
\]
Note that if $X=b_1E_1+\cdots+b_nE_n$ is a simple positive
operator in $\J(\M)$ and $|\!|\!|X|\!|\!|^\#\leq 1$, then
$\mu_s(X)\in \F'$ and $
|\!|\!|T|\!|\!|_{\mu_s(X)}=\int_0^\infty\mu_s(X)\mu_s(T)ds=c\sum_{k=1}^n
a_k^*b_k^*$, where $\{a_k^*\}$ and $\{b_k^*\}$ are nondecreasing
rearrangements of $\{a_k\}$ and $\{b_k\}$, respectively. By the
Hardy-Littlewood-Polya Theorem~\cite{HLP}, $\sum_{k=1}^n
a_kb_k\leq \sum_{k=1}^na_k^*b_k^*$. Hence,
\begin{multline*}
|\!|\!|T|\!|\!|=\sup\left\{c\sum_{k=1}^n a_kb_k:\,\,
X=b_1E_1+\cdots+b_nE_n\geq 0,\,|\!|\!|X|\!|\!|^\#\leq
1\right\}\\\leq\sup\{|\!|\!|T|\!|\!|_{f}:\, f\in
\F'\}=|\!|\!|T|\!|\!|'.
\end{multline*}

Now we need to  prove that $|\!|\!|T|\!|\!|'\leq |\!|\!|T|\!|\!|$.
Let $X\in \J(\M)$ be a positive simple operator such that
$|\!|\!|X|\!|\!|^\#\leq 1$. We need to prove that
$|\!|\!|T|\!|\!|_{\mu_s(X)}\leq |\!|\!|T|\!|\!|$. Since $T,X\in
\J(\M)$, there is a
finite projection $E\in \M$ such that $T,X\in \M_E$.\\

First, we assume that
$T=\tilde{a}_1\tilde{E}_1+\cdots+\tilde{a}_r\tilde{E}_r$ and
$X=b_1F_1+\cdots+b_rF_r$, where $E_1+\cdots+E_r=F_1+\cdots+F_r=E$,
$\tau(\tilde{E}_i)=\tau(F_j)=\tilde{c}$ for $1\leq i,j\leq r$,
$\tilde{a}_1\geq \cdots\geq \tilde{a}_r$ and $b_1\geq \cdots\geq
b_r$. Let $Y=b_1\tilde{E}_1+\cdots+b_r\tilde{E}_r$. Then
$\mu_s(Y)=\mu_s(X)$.  Since
$\tau(\tilde{E}_i)=\tau(F_j)=\tilde{c}$ for $1\leq i,j\leq r$ and
$\Aut(\M,\tau)$ acts on $\M$ ergodically, there is a $\theta\in
\Aut(\M,\tau)$ such that $\theta(\tilde{E_i})=F_i$ for $1\leq
i\leq r$. Hence $\theta(Y)=X$. Since $|\!|\!|\cdot|\!|\!|^\#$ is a
symmetric gauge norm, $|\!|\!|Y|\!|\!|^\#=|\!|\!|X|\!|\!|^\#\leq
1$. By Corollary~\ref{C:Holder inequality for bounded operators},
$|\!|\!|T|\!|\!|\geq\tau(TY)=\tilde{c}\sum_{k=1}^r
\tilde{a}_kb_k=\int_0^\infty\mu_s(Y)\mu_s(T)ds=\int_0^\infty\mu_s(X)\mu_s(T)ds=
|\!|\!|T|\!|\!|_{\mu_s(X)}$. \\

Now we consider the general case.
 Since $(\M_E,\tau_E)$ satisfies the weak Dixmier
property, by Lemma~\ref{L:weak dixmier property}, $\M_E$ is either
a finite dimensional von Neumann algebra such that
$\tau(F)=\tau(F')$ for every two minimal projections $F$ and $F'$
or $\M_E$ is a diffuse von Neumann algebra. The first case is
proved. If $\M_E$ is a diffuse von Neumann algebra, by
Lemma~\ref{L:proof of theorem c}, we can construct a sequence of
simple positive operators $X_n\in\J(\M)$ satisfying the following
conditions:
\begin{enumerate}
\item $0\leq X_1\leq X_2\leq \cdots\leq X$ and hence $0\leq
\mu_s(X_1)\leq \mu_s(X_2)\leq \cdots\leq \mu_s(X)$ for all $s\in
[0,\infty)$; \item $\lim_{n\rightarrow \infty}\mu_s(X_n)=\mu_s(X)$
for almost all $s\in [0,\infty)$; \item there exists an $r_n\in
\nn$ such that
$T={a}_{n,1}{E}_{n,1}+\cdots+{a}_{n,r_n}{E}_{n,r_n}$ and
$X=b_{n,1}F_{n,1}+\cdots+b_{n,r_n}F_{n,r_n}$, where
$E_{n,1}+\cdots+E_{n,r_n}=F_{n,1}+\cdots+F_{n,r_n}=E$ and
$\tau(E_{n,i})=\tau(F_{n,j})$ for $1\leq i,j\leq n$.
\end{enumerate}
By 1 and Corollary~\ref{C:S<T}, $|\!|\!|X_n|\!|\!|^\#\leq
|\!|\!|X|\!|\!|\leq 1$ for all $n=1,2,\cdots$. By 3 and the above
arguments, for every $n$, $|\!|\!|T|\!|\!|\geq
|\!|\!|T|\!|\!|_{\mu_{s}(X_n)}$. By 1, 2 and the Monotone
Convergence theorem,
\[|\!|\!|T|\!|\!|_{\mu_{s}(X)}=\int_0^\infty
\mu_s(X)\mu_s(T)ds=\lim_{n\rightarrow\infty}\int_0^\infty
\mu_s(X_n)\mu_s(T)ds=\lim_{n\rightarrow\infty}|\!|\!|T|\!|\!|_{\mu_{s}(X_n)}\leq
|\!|\!|T|\!|\!|.\]

 \end{proof}

\begin{Corollary}\label{C:extension}Let $(\M,\tau)$ be a semi-finite von Neumann
 algebra as in  Theorem~\ref{T:Theorem C} and let
 $|\!|\!|\cdot|\!|\!|$ be a normalized symmetric gauge norm on
$\J(\M)$. Then $|\!|\!|\cdot|\!|\!|$ can be extended to a normalized
symmetric gauge norm $|\!|\!|\cdot|\!|\!|'$ associated to $\M$.
\end{Corollary}
\begin{proof}For $T\in \M$, define $|\!|\!|T|\!|\!|'=\max\{|\!|\!|T|\!|\!|_f:\, f\in
\F'\}$. Then $|\!|\!|\cdot|\!|\!|'$ is an extension of
$|\!|\!|\cdot|\!|\!|$.
\end{proof}

\begin{Remark}\emph{ In Corollary~\ref{C:extension}, the extension is not
unique. Indeed, define $|\!|\!|\cdot|\!|\!|$ on $\B(\H)$ by
$|\!|\!|T|\!|\!|=\|T\|$ if $T$ is a finite rank operator and
$|\!|\!|T|\!|\!|=\infty$ if $T$ is an infinite rank operator. It is
easy to see that $|\!|\!|\cdot|\!|\!|$ defines a unitarily invariant
norm associated to $\B(\H)$ such that the restriction of
$|\!|\!|\cdot|\!|\!|$ on $\J(\H)$ is the operator norm.}
\end{Remark}

\section{Unitarily invariant norms related to semi-finite factors}

As the first application of Theorem~\ref{T:Theorem C}, we set up a
structure theorem for unitarily invariant norms related to
semi-finite factors. Recall that $\F$ is the set of
non-increasing, non-negative, right continuous simple functions
$f(x)$ on $[0,\infty)$  with compact supports  such that $\int_0^1
f(x)dx\leq 1$.
\begin{Theorem}\label{T:Theorem A} Let $\M$ be a semi-factor and
let $|\!|\!|\cdot|\!|\!|$ be a  unitarily invariant norm on
$\J(\M)$. Then there is a subset $\F'$ of $\F$ containing the
characteristic function on $[0,1]$ such that for all $T\in
\J(\M)$, $|\!|\!|T|\!|\!|=\sup\{|\!|\!|T|\!|\!|_{f}:\, f\in
\F'\}$.
\end{Theorem}
\begin{proof}Combining Theorem~\ref{T:unitarily
invariant norms on type III} and Proposition~\ref{P:unitarily
invariant implies trace-preserving}, we prove the theorem.
\end{proof}

The next corollary also follows from Theorem~\ref{T:Theorem C}.
\begin{Corollary}\label{C:representation theorem for UIN on 21 factor}
 Let $|\!|\!|\cdot|\!|\!|$ be a normalized symmetric gauge norm on
$\J(L^\infty[0,\infty))$. Then there is a subset $\F'$ of $\F$
containing the characteristic function on $[0,1]$ such that for
all $T\in \J(L^\infty[0,\infty))$,
\[|\!|\!|T|\!|\!|=\sup\{|\!|\!|T|\!|\!|_{f}:\, f\in \F'\}.\]
\end{Corollary}
\vskip 1cm

Let $\M=\B(\H)$ and $\tau={\rm Tr}$.  By the proof of
Theorem~\ref{T:Theorem C}, if $f\in \F'$, then $f(s)=\mu_s(X)$ for
some finite rank operator $X\in \B(\H)$, $X\geq 0$ and
$|\!|\!|X|\!|\!|^\#\leq 1$. Write $\mu_s(X)=s_1(X)
\chi_{[0,1)}(s)+s_2(X)\chi_{[1, 2)}(s)+\cdots$, where $s_1(X),
s_2(X),\cdots,$ are  $s$-numbers of $X$. Since $\int_0^1
\mu_s(X)\leq 1$, $s_1(X)\leq 1$. By Lemma~\ref{L:T_f} and simple
computations, for every $T\in \J(\H)$,
\[|\!|\!|T|\!|\!|_{\mu_s(X)}=s_1(X) s_1(T)+s_2(X)
s_2(T)+\cdots,\] where $s_1(T),s_2(T),\cdots,$ are  $s$-numbers of $T$.\\

Let $\G=\{(a_1,a_2,\cdots):\, 1\geq a_1\geq a_2\geq \cdots\geq
0\,\text{and $a_n=0$ except for finitely many terms}\}$. For
$(a_1,a_2,\cdots)\in \G$ and $T\in \J(\H)$, define
\begin{equation}\label{Eq: a1 a2}
 |\!|\!|T|\!|\!|_{(a_1,\,a_2,\cdots)}=a_1s_1(T)+a_2s_2(T)+\cdots. \end{equation} Then
$|\!|\!|T|\!|\!|_{(a_1,a_2,\cdots)}=|\!|\!|T|\!|\!|_{f}$ is a
unitarily invariant norm on $\J(\H)$, where $f(x)=a_1
\chi_{[0,1)}(x)+a_2\chi_{[1, 2)}(x)+\cdots$. By identifying
$\mu_s(X)$ with $(s_1(X),s_2(X),\cdots)$ in $\G$, we obtain the
following corollary.
\begin{Corollary}\label{C:representation theorem for MGN on finite dimensional}
 Let  $|\!|\!|\cdot|\!|\!|$
be a unitarily invariant norm on $\J(\H)$. Then there is a subset
$\G'$ of $\G$, $(1,0,\cdots)\in \G'$, such that for all $T\in
\J(\H)$,
\[ |\!|\!|T|\!|\!|=\sup\{a_1s_1(T)+a_2s_2(T)+\cdots:\, (a_1, a_2, \cdots)\in
\G'\},
\] where $s_1(T), s_2(T), \cdots$ are $s$-numbers of $T$.
\end{Corollary}

Similar to the proof of Corollary~\ref{C:representation theorem for
MGN on finite dimensional}, we have the following corollary.
\begin{Corollary}\label{C:representation theorem for UIN on Cn}
 Let $|\!|\!|\cdot|\!|\!|$ be a normalized symmetric gauge
norm on $c_{00}=\J(l^\infty(\nn))$. Then there is a subset $\G'$
of $\G$, $(1,0,\cdots)\in \G'$, such that for all
$(x_1,x_2,\cdots)\in c_{00}$,
\[ |\!|\!|(x_1,x_2,\cdots)|\!|\!|=\sup\{a_1x_1^*+a_2x_2^*+\cdots:\, (a_1,a_2,\cdots)\in
\G'\},
\] where $(x_1^*,x_2^*,\cdots)$ is the nonincreasing rearrangement of
$(|x_1|, |x_2|,\cdots)$.
\end{Corollary}

\section{Unitarily invariant norms and symmetric gauge norms}

\begin{Lemma}\label{L: two embedding are equivalent}
 Let $\theta_1, \theta_2$ be two embeddings from
$(L^\infty[0,\infty), \int_0^\infty dx)$ into a type ${\rm II}\sb
\infty$ factor $(\M,\tau)$. If $|\!|\!|\cdot|\!|\!|$ is a unitarily
invariant norm on $\J(\M)$, then
$|\!|\!|\theta_1(f)|\!|\!|=|\!|\!|\theta_2(f)|\!|\!|$ for every
positive function $f\in \J(L^\infty[0,\infty))$.
\end{Lemma}
\begin{proof}  For $f\in \J(L^\infty[0,\infty))$, let $|\!|\!|f|\!|\!|_1=|\!|\!|\theta_1(f)|\!|\!|$ and
$|\!|\!|f|\!|\!|_2=|\!|\!|\theta_2(f)|\!|\!|$.  Then
$|\!|\!|\cdot|\!|\!|_1$ and $|\!|\!|\cdot|\!|\!|_2$ are gauge norms
on $\J(L^\infty[0,\infty))$. By Lemma~\ref{L:simple operators are
dense}, to prove $|\!|\!|\cdot|\!|\!|_1=|\!|\!|\cdot|\!|\!|_2$ on
$\J(L^\infty[0,\infty))$, we need to prove
$|\!|\!|f|\!|\!|_1=|\!|\!|f|\!|\!|_2$ for every simple function
$f(x)$ in $\J(L^\infty[0,\infty))$. If $f(x)\in
\J(L^\infty[0,\infty))$ is a simple function, then there is a
unitary operator $U$ in $\M$ such that
$U\theta_1(f)U^*=\theta_2(f)$. Hence
$|\!|\!|f|\!|\!|_1=|\!|\!|f|\!|\!|_2$.
\end{proof}

The following theorem generalizes von Neumann's classical
result~\cite{vN} on unitarily invariant norms on $M_n(\cc)$.

\begin{Theorem}\label{T:Theorem D}
 There is a one-to-one correspondence between
\begin{enumerate}
\item unitarily invariant norms on $M_n(\cc)$ and symmetric gauge
norms on $\cc^n$, \item unitarily invariant norms on a type ${\rm
II}_1$ factor and symmetric gauge norms on $L^\infty[0,1]$,\item
unitarily invariant norms on $\J(\H)$ and symmetric gauge norms on
$c_{00}=\J(l^\infty(\mathbb{N}))$,\item unitarily invariant norms
on $\J(\M)$ of a type ${\rm II}_\infty$ factor $\M$ and symmetric
gauge norms on $\J(L^\infty[0,\infty))$,
\end{enumerate}
 respectively. Precisely, let
$\M$ be a semi-finite factor and $\A$ be the corresponding abelian
von Neumann algebra as the above:
\begin{itemize}
\item if $|\!|\!|\cdot|\!|\!|$ is a unitarily invariant norm on
$\J(\M)$ and $\theta$ is an embedding from  $\A$ into  $\M$, then
the restriction of $|\!|\!|\cdot|\!|\!|$ to $\J(\theta(\A))$
defines a symmetric gauge norm on $\J(\A)$; \item conversely, if
$|\!|\!|\cdot|\!|\!|'$ is a symmetric gauge norm on $\J(\A)$ and
$T\in \J(\M)$, then $|\!|\!|\mu_s(T)|\!|\!|'$ defines a unitarily
invariant norm on $\J(\M)$, where $\mu_s(T)$ is the classical
$s$-numbers of $T$ if $\M=M_n(\cc)$ or $\M=\B(\H)$, and $\mu_s(T)$
is defined as in~\cite{FHNS} if $\M$ is a type ${\rm II}_1$
factor.
\end{itemize}
\end{Theorem}
\begin{proof} We refer to~\cite{FHNS} for the proof of case 1 and
case 2.
 We only prove the
theorem for case 4. The proof of case 3 is similar.  We may assume
that both norms on $\J(\M)$ and $\J(L^\infty[0,\infty))$ are
normalized. By the definition of Ky Fan norms, there is a
one-to-one correspondence
 between  Ky Fan $t$-th norms on $\J(\M)$ and  Ky Fan $t$-th norms on
$\J(L^\infty[0,\infty))$ as in Theorem~\ref{T:Theorem D}. By
Theorem~\ref{T:Theorem A} and Corollary~\ref{C:representation
theorem for UIN on 21 factor}, there is a one-to-one
correspondence between normalized unitarily invariant norms on
$\J(\M)$ and normalized symmetric gauge norms on
$\J(L^\infty[0,\infty))$ as in the theorem.
\end{proof}

\begin{Example}\label{E:Lp space} \emph{For $1\leq p\leq \infty$, the $L^p$-norm  on
$(L^\infty[0,\infty), \int_0^\infty dx)$ defined by
\[\|f(x)\|_p=\left\{
               \begin{array}{ll}
                 \left(\int_0^\infty|f(x)|^pdx\right)^{1/p}, & \hbox{$1\leq p<\infty$;} \\
                 \text{essup}\, |f|, & \hbox{$p=\infty$}
               \end{array}
             \right.
\] is a normalized symmetric gauge norm on $(L^\infty[0,\infty), \int_0^\infty\,dx)$.
 By Theorem~\ref{T:Theorem D}, the
induced norm for $T\in \J(\M)$ of a type ${\rm II}\sb \infty$ factor
$\M$ defined by
\[\|T\|_p=\left\{
               \begin{array}{ll}
                =\left(\tau(|T|^p)\right)^{1/p}= \left(\int_0^1|\mu_s(T)|^pds\right)^{1/p}, & \hbox{$1\leq p<\infty$;} \\
                 \|T\|, & \hbox{$p=\infty$}
               \end{array}
             \right.
\] is a unitarily invariant norm on  $\J(\M)$.
The norms $\{\|\cdot\|_p:\, 1\leq p\leq \infty\}$ are called
$L^p$-norms on $\J(\M)$.}
\end{Example}
\begin{Example}\label{E:Lp for discrete} \emph{For $1\leq p\leq \infty$, the $l^p$-norm defined on
$\J(l^\infty(\nn))$ by
\[\|(x_1,x_2,\cdots)\|_p=\left\{
               \begin{array}{ll}
                 \left(|x_1|^p+|x_2|^p+\cdots\right)^{1/p}, & \hbox{$1\leq p<\infty$;} \\
                 \text{essup}\, \{|x_n|:\, n=1,2,\cdots\}, & \hbox{$p=\infty$}
               \end{array}
             \right.
\] is a normalized symmetric gauge norm on $\J(l^\infty(\nn))$. By Theorem~\ref{T:Theorem
D}, the induced norm for $T$ in $\J(\H)$ defined by
\[\|T\|_p=\left\{
               \begin{array}{ll}
                =\left(\tau(|T|^p)\right)^{1/p}= \left(s_1(T)^p+s_2(T)^p+\cdots\right)^{1/p}, & \hbox{$1\leq p<\infty$;} \\
                 \|T\|, & \hbox{$p=\infty$}
               \end{array}
             \right.
\] is a unitarily invariant norm on $\J(\H)$. The norms
$\{\|\cdot\|_p:\, 1\leq p\leq \infty\}$ are called $L^p$-norms on
$\J(\H)$.}
\end{Example}

Theorem~\ref{T:Theorem D}  establishes  the one to one
correspondence between unitarily invariant norms on $\J(\M)$ of an
infinite semi-finite factor $\M$ and symmetric gauge norms on
$\J(\A)$ of an abelian von Neumann algebra $\A$. The following
theorem further establishes the one to one correspondence between
the dual norms on $\J(\M)$ and the dual norms on $\J(\A)$, which
plays a key role in the studying of duality and reflexivity of the
completion  of $\J(\M)$ with respect to  unitarily invariant
norms.
\begin{Theorem}\label{T:Theorem E}
Let $\M$ be a type ${\rm II}\sb \infty$  factor $($or $\B(\H)$$)$.
If $|\!|\!|\cdot|\!|\!|$ is a unitarily invariant norm on $\J(\M)$
corresponding to the symmetric gauge norm $|\!|\!|\cdot|\!|\!|_1$
on $\J(L^\infty[0,\infty))$ $($or $\J(l^\infty(\nn))$
respectively$)$ as in Theorem~\ref{T:Theorem D}, then
$|\!|\!|\cdot|\!|\!|^\#$ is the unitarily invariant norm  on
$\J(\M)$ corresponding to the symmetric gauge norm
$|\!|\!|\cdot|\!|\!|_1^\#$ on $\J(L^\infty[0,\infty))$ $($or
$\J(l^\infty(\nn))$ respectively$)$ as in Theorem~\ref{T:Theorem
D}.
\end{Theorem}
\begin{proof} We only prove the theorem for type
${\rm II}\sb \infty$ factors. The case of type ${\rm I}\sb \infty$
factors is similar. Let
 $|\!|\!|\cdot|\!|\!|_2$ be the unitarily invariant norm on
 $\J(\M)$ corresponding to the symmetric gauge norm  $|\!|\!|\cdot|\!|\!|_1^\#$
on $\J(L^\infty[0,\infty))$  as in Theorem~\ref{T:Theorem D}. By
Lemma~\ref{L:simple operators are dense}, to prove
$|\!|\!|\cdot|\!|\!|_2=|\!|\!|\cdot|\!|\!|^\#$ on $\J(\M)$, we
need to prove $|\!|\!|T|\!|\!|_2=|\!|\!|T|\!|\!|^\#$ for every
simple positive operator $T=a_1E_1+\cdots+a_nE_n$ in $\J(\M)$ such
that $\tau(E_1)=\cdots=\tau(E_n)=c$. We may assume that $a_1\geq
\cdots\geq a_n\geq 0$. Then
$\mu_s(T)=a_1\chi_{[0,c)}(s)+\cdots+a_n\chi_{[(n-1)c,nc)}(s)$. By
Lemma~\ref{L:Dual norm of simple operators},
\[|\!|\!|T|\!|\!|^\#=\sup\{c\sum_{k=1}^na_kb_k:\,
X=b_1E_1+\cdots+b_nE_n\geq 0,\, |\!|\!|X|\!|\!|\leq 1\}.\] By the
Hardy-Littlewood-Polya Theorem,
\[|\!|\!|T|\!|\!|^\#=\sup\{c\sum_{k=1}^na_kb_k:\,
X=b_1E_1+\cdots+b_nE_n\geq 0,\, b_1\geq \cdots\geq b_n\geq 0,\,
|\!|\!|X|\!|\!|\leq 1\}.\] By Theorem~\ref{T:Theorem D} and
Lemma~\ref{L:Dual norm of simple operators},
\[|\!|\!|T|\!|\!|_2=|\!|\!|\mu_s(T)|\!|\!|^\#=\sup\{c\sum_{k=1}^na_kb_k:\,
g(s)=b_1\chi_{[0,c)}(s)+\cdots+b_n\chi_{[(n-1)c,nc)}(s)\geq 0,\,
|\!|\!|g(s)|\!|\!|\leq 1\}.\] By the Hardy-Littlewood-Polya Theorem,
\begin{multline*}
|\!|\!|T|\!|\!|_2=|\!|\!|\mu_s(T)|\!|\!|^\#=\sup\{c\sum_{k=1}^na_kb_k:\,
g(s)=b_1\chi_{[0,c)}(s)+\cdots+b_n\chi_{[(n-1)c,nc)}(s)\geq 0,\\
b_1\geq \cdots b_n\geq 0,\, |\!|\!|g(s)|\!|\!|\leq 1\}.
\end{multline*}
 Note that if
$b_1\geq \cdots\geq b_n\geq 0$, then
$\mu_s(b_1E_1+\cdots+b_nE_n)=b_1\chi_{[0,c)}(s)+\cdots+b_n\chi_{[(n-1)c,nc)}(s)$.
Since $|\!|\!|\cdot|\!|\!|$ is the unitarily invariant norm on
$\J(\M)$ corresponding to the symmetric gauge norm
$|\!|\!|\cdot|\!|\!|_1$ on $\J(L^\infty[0,\infty))$ as in
Theorem~\ref{T:Theorem D}, $|\!|\!|b_1E_1+\cdots+b_nE_n|\!|\!|\leq
1$ if and only if
$|\!|\!|b_1\chi_{[0,c)}(s)+\cdots+b_n\chi_{[(n-1)c,nc)}(s)|\!|\!|_1\leq
1$. Therefore, $|\!|\!|T|\!|\!|_2=|\!|\!|T|\!|\!|^\#$.
\end{proof}

\begin{Example}\label{E:dual Lp 1}\emph{ If $p=1$, let $q=\infty$. If $1<p<\infty$, let $q=\frac{p}{p-1}$.
Then the $L^q$ norm  on $\J(L^\infty[0,\infty))$ defined by
Example~\ref{E:Lp space} is the dual norm of the $L^p$-norm on
$\J(L^\infty[0,\infty))$. By Theorem~\ref{T:Theorem E}, the
$L^q$-norm on $\J(\M)$ of a type ${\rm II}\sb \infty$ factor $\M$
is the dual norm of the $L^p$-norm on $\J(\M)$.}
\end{Example}

\begin{Example}\label{E:dual Lp 2}\emph{ If $p=1$, let $q=\infty$. If $1<p<\infty$, let $q=\frac{p}{p-1}$.
Then the $l^q$-norm  on $\J(l^\infty(\nn))$ defined by
Example~\ref{E:Lp for discrete} is the dual norm of the $l^p$-norm
on $\J(l^\infty(\nn))$. By Theorem~\ref{T:Theorem E}, the
$L^q$-norm on $\J(\H)$ is the dual norm of the $L^p$-norm on
$\J(\H)$.}
\end{Example}

\section{Ky Fan's dominance theorem}

The following theorem generalizes Ky Fan's dominance theorem.
\begin{Theorem}\label{T:Theorem F}
Let $\M$ be a semi-finite factor
 and $S, T\in \J(\M)$. If
$|\!|\!|S|\!|\!|_{(t)}\leq |\!|\!|T|\!|\!|_{(t)}$ for all Ky Fan
$t$-th norms, $0\leq t\leq \infty$, then $|\!|\!|S|\!|\!|\leq
|\!|\!|T|\!|\!|$ for all unitarily invariant norms
$|\!|\!|\cdot|\!|\!|$ on $\J(\M)$.
\end{Theorem}

\begin{proof}  Let $|\!|\!|\cdot|\!|\!|$ be a unitarily invariant norm on
$\M$. By Lemma~\ref{L:T_f}, $|\!|\!|S|\!|\!|_f\leq
|\!|\!|T|\!|\!|_f$ for every $f\in \F$. By Theorem~\ref{T:Theorem
A}, $|\!|\!|S|\!|\!|\leq |\!|\!|T|\!|\!|$.
\end{proof}

\begin{Corollary}\label{T:Ky Fan Theorem} Let  $S, T\in \J(\H)$. If
$|\!|\!|S|\!|\!|_{(n)}\leq |\!|\!|T|\!|\!|_{(n)}$ for all Ky Fan
$n$-th norms, $n=1,2,\cdots$, then $|\!|\!|S|\!|\!|\leq
|\!|\!|T|\!|\!|$ for all unitarily invariant norms
$|\!|\!|\cdot|\!|\!|$ on $\J(\H)$.
\end{Corollary}

\section{Completion of $\J(\M)$ with respect to unitarily invariant norms}
In this section, we assume that $\M$ is an infinite semi-finite
factor with a tracial weight $\tau$ and $|\!|\!|\cdot|\!|\!|$ is a
unitarily invariant norm on $\J(\M)$. Recall that $\J(\M)$ is the
set of operators $T$ in $\M$ such that $T=ETE$ for some finite
projection $E\in \M$.
 The completion of $\J(\M)$ with
respect to $|\!|\!|\cdot|\!|\!|$ is denoted by
$\overline{{\J(\M)}_{|\!|\!|\cdot|\!|\!|}}$.  For $1\leq p<\infty$,
we will use the traditional notation $L^p(\M,\tau)$ to denote the
completion of $\J(\M)$ with respect to the $L^p$-norm defined as in
Example~\ref{E:Lp space} or Example~\ref{E:Lp for discrete}. We will
denote by $\K(\M)$ the completion of $\J(\M)$ with respect to the
operator norm.  Let $\widetilde{\M}$ be the completion of $\M$ in
the measure-topology in the sense of Nelson~\cite{Ne}. Recall that a
neighborhood $N(\epsilon,\delta)$ of $0\in\M$ in the measure
topology (see~\cite{Ne}) is
\[N(\epsilon,\delta)=\{T\in \M, \,\text{there is a projection $E\in\M$ such that
$\tau(E)<\delta$ and $\|TE^\perp\|<\epsilon$}\}
\]

\subsection{$L^1(\M,\tau)$}

In this section, we summarize some classical results on
$\widetilde{\M}$ and $L^1(\M,\tau)$ that will be useful.
\begin{Theorem}\label{T:Nelson}\emph{(Nelson, \cite{Ne})} $\widetilde{\M}$
is a $\ast$-algebra and every element in $\widetilde{\M}$ is a
closed, densely defined operator affiliated with $\M$.
\end{Theorem}

In the following, we  define $s$-numbers for unbounded operators in
$\widetilde{\M}$ as~\cite{F-K}.
\begin{Definition}\emph{ For $T\in \widetilde{\M}$ and $0\leq s<\infty$, define
the \emph{$s$-numbers of $T$} by
\[\mu_s(T)=\inf\{\|TE\|:\,\,\text{$E\in \M$ is a projection such that $\tau(E^\perp)=s$
}\}.
\]}
\end{Definition}

\begin{Theorem}\emph{(Fack and   Kosaki, \cite{F-K})}
\label{T:THIERRY FACK AND HIDEKI KOSAKI} Let $T$ and $T_n$ be a
sequence of operators in $\widetilde{\M}$ such that
$\lim_{n\rightarrow \infty}T_n=T$ in the measure-topology. Then
for almost all $s\in [0,\infty)$,
$\lim_{n\rightarrow\infty}\mu_s(T_n)=\mu_s(T)$.
\end{Theorem}

Let $\{T_n\}$ be a sequence of operators in $\J(\M)$ such that
$T=\lim_{n\rightarrow\infty} T_n$ in the $L^1$-norm. By
Lemma~\ref{L:characterization of Ky fan norms 2}, $\{\tau(T_n)\}$ is
a Cauchy sequence in $\cc$. Define
$\tau(T)=\lim_{n\rightarrow\infty}\tau(T_n)$. It is obvious that
$\tau(T)$ does not depend on the sequence $\{T_n\}$. In this way,
$\tau$ is extended to a linear functional on $L^1(\M,\tau)$. The
following lemma is due to Nelson~\cite{Ne}.
\begin{Lemma}\label{C:embedding L1 into measure topology} Let
$\M$ be a type ${\rm II}\sb \infty$ factor with a faithful normal
tracial weight $\tau$. Then there is a natural continuous injective
map from $L^1(\M,\tau)$ into $\widetilde{\M}$. Furthermore, if
$\{T_n\}\subset\J(\M)$ is a Cauchy sequence in the $L^1$-norm and
$\lim_{n\rightarrow\infty}T_n=T$ in the measure topology, then $T\in
L^1(\M,\tau)$ and $\lim_{n\rightarrow\infty}T_n=T$ in the
$L^1$-norm.
\end{Lemma}

\subsection{Embedding of $\overline{{\J(\M)}_{|\!|\!|\cdot|\!|\!|}}$ into $\widetilde{\M}$}

In this subsection, we will prove that
$\overline{{\J(\M)}_{|\!|\!|\cdot|\!|\!|}}$ can be embedded into
$\widetilde{\M}$. If $\M=\B(\H)$, the embedding is very simple due
to the following proposition.
\begin{Proposition}\label{P:closure of J(H)}Let $|\!|\!|\cdot|\!|\!|$ be a unitarily invariant norm on
$\J(\H)$. Then $\overline{{\J(\H)}_{|\!|\!|\cdot|\!|\!|}}$ is an
selfadjoint two-sided ideal of $\B(\H)$ such that $\K_1(\H)\subseteq
\overline{{\J(\H)}_{|\!|\!|\cdot|\!|\!|}}\subseteq \K(\H)$.
\end{Proposition}
\begin{proof} By Corollary~\ref{C:comparison corollary}, we obtain
the proposition.
\end{proof}

\begin{Lemma}\label{L:embedding into measure topology} Let
$\M$ be a type ${\rm II}\sb \infty$ factor and $|\!|\!|\cdot|\!|\!|$
be a unitarily invariant norm on $\J(\M)$. There is a natural
continuous map $\Phi$ from
$\overline{{\J(\M)}_{|\!|\!|\cdot|\!|\!|}}$ to $\widetilde{\M}$ that
extends the identity map from $\J(\M)$ to $\J(\M)$.
\end{Lemma}
\begin{proof}  We may assume that $|\!|\!|\cdot|\!|\!|$ is a
normalized unitarily invariant norm, i.e., $|\!|\!|E|\!|\!|=1$ if
$\tau(E)=1$.
 If
$\{T_n\}$ in $\J(\M)$ is a Cauchy sequence with respect to
$|\!|\!|\cdot|\!|\!|$, then $\{T_n\}$ in $\J(\M)$ is a Cauchy
sequence with respect to $|\!|\!|\cdot|\!|\!|_{(1)}$ by
Proposition~\ref{P:comparison proposition}.  By
Corollary~\ref{C:Markov inequality}, for any $\delta>0$ and $T\in
\J(\M)$ such that $|\!|\!|T|\!|\!|_{(1)}<\delta$,
$\tau(\chi_{(\delta,\infty)}(|T|))\leq
\frac{|\!|\!|T|\!|\!|_{(1)}}{\delta}$. Hence, if $\{T_n\}$ is a
Cauchy sequence in $\M$ with respect to the
$|\!|\!|T|\!|\!|_{(1)}$ norm, then $\{T_n\}$ is a Cauchy sequence
in the measure topology. So there is a natural map $\Phi$ from
$\overline{{\J(\M)}_{|\!|\!|\cdot|\!|\!|}}$ that extends the
identity map from $\J(\M)$ to $\J(\M)$.
\end{proof}

\begin{Lemma}\label{L:convergence in L1 norm}Let $\{T_n\}$ be a sequence of operators in $\J(\M)$ such
that $\lim_{n\rightarrow\infty} T_n=T$ with respect to
$|\!|\!|\cdot|\!|\!|$ and $X\in \J(\M)$. Then
$\lim_{n\rightarrow\infty}T_nX=\Phi(T)X$ in the $L^1$-norm.
\end{Lemma}
\begin{proof}By Lemma~\ref{L:embedding into measure topology},
$\lim_{n\rightarrow\infty}T_n=\Phi(T)$ in the measure topology.
Hence $\lim_{n\rightarrow\infty}T_nX=\Phi(T)X$ in the measure
topology (see Theorem 1 of~\cite{Ne}). By Corollary~\ref{C:Holder
inequality for bounded operators}, $\|T_nX-T_mX\|_1\leq
|\!|\!|T_n-T_m|\!|\!|\cdot |\!|\!|X|\!|\!|^\#$. Hence $\{T_nX\}$ is
a Cauchy sequence in the $L^1$-norm. By Lemma~\ref{C:embedding L1
into measure topology}, $\lim_{n\rightarrow\infty}T_nX=\Phi(T)X$ in
the $L^1$-norm.
\end{proof}
\begin{Corollary}\label{C:convergence in L1 norm}Let $\{T_n\}$ be a sequence of operators in $\J(\M)$ such
that $\lim_{n\rightarrow\infty} T_n=T$ with respect to
$|\!|\!|\cdot|\!|\!|$ and $X\in \J(\M)$. Then
$\lim_{n\rightarrow\infty}\tau((T_n-\Phi(T))X)=0$.
\end{Corollary}
\begin{Lemma}\label{L:dual norms for unbounded operators} For all $T\in \overline{{\J(\M)}_{|\!|\!|\cdot|\!|\!|}}$,
\[|\!|\!|T|\!|\!|=\sup\{|\tau(\Phi(T)X)|:\, X\in\J(\M),\, |\!|\!|X|\!|\!|^\#\leq 1\}.
\]
\end{Lemma}
\begin{proof}Let $\{T_n\}$ be a sequence of operators in $\J(\M)$ such
that $\lim_{n\rightarrow\infty} T_n=T$ with respect to
$|\!|\!|\cdot|\!|\!|$. By Corollary~\ref{C:convergence in L1 norm}
and Lemma~\ref{L:Holder inequality for bounded operators}, for every
$X\in \J(\M)$ such that $|\!|\!|X|\!|\!|^\#\leq 1$,
\[|\tau(\Phi(T)X)|=\lim_{n\rightarrow \infty}|\tau(T_nX)|\leq
\lim_{n\rightarrow \infty} |\!|\!|T_n|\!|\!|=|\!|\!|T|\!|\!|.\]
Hence, $|\!|\!|T|\!|\!|\geq \sup\{|\tau(\Phi(T)X)|:\, X\in\J(\M),\,
|\!|\!|X|\!|\!|^\#\leq 1\}.$\\

 We need to prove that $|\!|\!|T|\!|\!|\leq \sup\{|\tau(\Phi(T)X)|:\, X\in\J(\M),\,
|\!|\!|X|\!|\!|^\#\leq 1\}.$ Let $\epsilon>0$. Since
$\lim_{n\rightarrow\infty} T_n=T$ with respect to
$|\!|\!|\cdot|\!|\!|$, there exists an $N$ such that
 $|\!|\!|T-T_N|\!|\!|<\epsilon/3$.
 For $T_N$, there is an
$X\in\J(\M)$, $|\!|\!|X|\!|\!|^\#\leq 1$, such that
$|\!|\!|T_N|\!|\!|\leq |\tau(T_NX)|+\epsilon/3$. By
Corollary~\ref{C:convergence in L1 norm} and Corollary~\ref{C:Holder
inequality for bounded operators}, \[|\tau((T_N-\Phi(T))X)|
=\lim_{n\rightarrow\infty}|\tau((T_N-T_n)X)\leq
\lim_{n\rightarrow\infty}|\!|\!|T_N-T_n|\!|\!|\cdot
|\!|\!|X|\!|\!|^\#\leq |\!|\!|T_N-T|\!|\!| <\epsilon/3.\]
 So
$|\tau(\Phi(T)X)|\geq |\tau(T_NX)|-|\tau((T_N-\Phi(T))X)|\geq
|\!|\!|T_N|\!|\!|-\epsilon/3-\epsilon/3\geq
|\!|\!|T|\!|\!|-\epsilon$. Therefore, $|\!|\!|T|\!|\!|\leq
\sup\{|\tau(\Phi(T)X)|:\, X\in\J(\M),\, |\!|\!|X|\!|\!|^\#\leq 1\}.$
\end{proof}

\begin{Proposition}\label{P:embedding into measure topology} Let
$\M$ be a type ${\rm II}\sb \infty$ factor and $|\!|\!|\cdot|\!|\!|$
be a unitarily invariant norm on $\J(\M)$. There is a natural
continuous injective map $\Phi$ from
$\overline{{\J(\M)}_{|\!|\!|\cdot|\!|\!|}}$ to $\widetilde{\M}$ that
extends the identity map from $\J(\M)$ to $\J(\M)$.
\end{Proposition}
\begin{proof} By Lemma~\ref{L:embedding into measure topology}, we
need only to prove that $\Phi$ is injective. To prove $\Phi$ is
injective, we need to prove that if $\{T_n\}$ in $\J(\M)$ is a
Cauchy sequence with respect to $|\!|\!|\cdot|\!|\!|$ and
$T_n\rightarrow 0$ in the measure topology, then
$\lim_{n\rightarrow\infty}|\!|\!|T_n|\!|\!|=0$. Suppose
$T=\lim_{n\rightarrow\infty}T_n$ with respect to the norm
$|\!|\!|\cdot|\!|\!|$. Then $\Phi(T)=0$. By Lemma~\ref{L:dual norms
for unbounded operators},
\[|\!|\!|T|\!|\!|=\sup\{|\tau(\Phi(T)X)|:\, X\in\J(\M),\, |\!|\!|X|\!|\!|^\#\leq 1\}=0.
\] Hence, $T=0$.
\end{proof}

By Proposition~\ref{P:embedding into measure topology}, we will
identify $T\in \overline{{\J(\M)}_{|\!|\!|\cdot|\!|\!|}}$ with
$\Phi(T)$ and consider $\overline{{\J(\M)}_{|\!|\!|\cdot|\!|\!|}}$
as a subset of $\widetilde{\M}$.

\begin{Corollary}\label{C:extension of norm on the completion algebra}
 $\overline{{\J(\M)}_{|\!|\!|\cdot|\!|\!|}}$ is a
 linear subspace of $\widetilde{\M}$ satisfying the following conditions:
\begin{enumerate}
\item if $T\in \overline{{\J(\M)}_{|\!|\!|\cdot|\!|\!|}}$, then
$T^*\in \overline{{\J(\M)}_{|\!|\!|\cdot|\!|\!|}}$;
\item $T\in \overline{{\J(\M)}_{|\!|\!|\cdot|\!|\!|}}$ if and only if
$|T|\in \overline{{\J(\M)}_{|\!|\!|\cdot|\!|\!|}}$;
\item if $T\in \overline{{\J(\M)}_{|\!|\!|\cdot|\!|\!|}}$ and $A, B\in
\J(\M)$,
 then $ATB\in \overline{{\J(\M)}_{|\!|\!|\cdot|\!|\!|}}$ and
$|\!|\!|ATB|\!|\!|\leq \|A\|\cdot |\!|\!|T|\!|\!|\cdot\|B\|$.
\end{enumerate}
In particular, $|\!|\!|\cdot|\!|\!|$ can be extended to a unitarily
invariant norm, also denoted by $|\!|\!|\cdot|\!|\!|$, on
$\overline{{\J(\M)}_{|\!|\!|\cdot|\!|\!|}}$.
\end{Corollary}

\subsection{Elements in $\overline{{\J(\M)}_{|\!|\!|\cdot|\!|\!|}}$}

The following theorem generalizes Theorem~\ref{T:Theorem C}. Its
proof   is based on Lemma~\ref{L:dual norms for unbounded
operators} and is similar to the proof of Theorem~\ref{T:Theorem
C}. So we omit the proof. Recall that $\F$ is the set of
non-increasing, non-negative, right continuous simple functions
$f(x)$ on $[0,\infty)$  with compact supports  such that $\int_0^1
f(x)dx\leq 1$.

\begin{Theorem}\label{T: a representation theorem 2}  If
 $|\!|\!|\cdot|\!|\!|$ is a normalized unitarily invariant norm on an infinite
 semi-finite  factor
$\J(\M)$, then there is a subset $\F'$ of $\F$ containing the
characteristic function on $[0,1]$ such that for all $T\in
\overline{{\J(\M)}_{|\!|\!|\cdot|\!|\!|}}$,
\[|\!|\!|T|\!|\!|=\sup\{|\!|\!|T|\!|\!|_{f}:\, f\in \F'\}.\]
\end{Theorem}

Combining Theorem~\ref{T: a representation theorem 2} and
Theorem~\ref{T:THIERRY FACK AND HIDEKI KOSAKI}, we have the
following corollaries.
\begin{Corollary}Let
 $|\!|\!|\cdot|\!|\!|$ be a unitarily invariant norm on $\J(\H)$ and
 $|\!|\!|\cdot|\!|\!|'$ be the corresponding symmetric gauge norm on
 $c_{00}=\J(l^\infty(\nn,\tau))$. If $T\in \K(\H)$, then $T\in
 \overline{{\J(\H)}_{|\!|\!|\cdot|\!|\!|}}$ if and only if $(\mu_1(T),\mu_2(T),\cdots)$ is
 in the completion of $c_{00}$ with respect to
 $|\!|\!|\cdot|\!|\!|'$. In this case,
 $|\!|\!|T|\!|\!|=|\!|\!|\mu_s(T)|\!|\!|'$.
\end{Corollary}

\begin{Corollary}Let $\M$ be a type ${\rm II}\sb \infty$ factor and
 $|\!|\!|\cdot|\!|\!|$ be a unitarily invariant norm on
$\J(\M)$ and
 $|\!|\!|\cdot|\!|\!|'$ be the corresponding symmetric gauge norm on
 $\J(L^\infty[0, \infty))$. If $T\in \widetilde{\M}$, then $T\in
 \overline{{\J(\M)}_{|\!|\!|\cdot|\!|\!|}}$ if and only if
 $\mu_s(T)\in \overline{{\J(L^\infty[0,\infty))}_{|\!|\!|\cdot|\!|\!|'}}$. In this case,
 $|\!|\!|T|\!|\!|=|\!|\!|\mu_s(T)|\!|\!|'$.
\end{Corollary}

\begin{Example}\emph{Let $T\in \K(\H)$ and $1\leq p\leq \infty$. Then
$T\in L^p(\H,Tr)$ if and only if $(\mu_1(T),\mu_2(T),\cdots)\in
l^p(\nn,\tau)$. In this case,
$\|T\|_p=\left(s_1(T)^p+s_2(T)^p+\cdots\right)^{1/p}$.}
\end{Example}

\begin{Example}\emph{Let $T\in \widetilde{\M}$ and $1\leq p\leq \infty$. Then
$T\in L^p(\M,\tau)$ if and only if $\mu_s(T)\in L^p[0,\infty)$. In
this case,
$\|T\|_p=\left(\int_0^1\mu_s(T)^pds\right)^{1/p}=\left(\int_0^\infty
\lambda^p\,d\mu_{|T|}\right)^{1/p}$.}
\end{Example}

\subsection{A generalization of H$\ddot{\text{o}}$lder's inequality}
The following theorem is a generalization of
H$\ddot{\text{o}}$lder's inequality.
\begin{Theorem}\label{T:Holder inequality} Let $\M$ be an infinite
semi-finite factor and $|\!|\!|\cdot|\!|\!|$ be a normalized
unitarily invariant norm on  $\J(\M)$. If
 $T\in \overline{{\J(\M)}_{|\!|\!|\cdot|\!|\!|}}$ and $S\in
 \overline{{\J(\M)}_{|\!|\!|\cdot|\!|\!|^\#}}$,
  then $TS\in L^1(\M,\tau)$ and $\|TS\|_1\leq
|\!|\!|T|\!|\!| \cdot |\!|\!|S|\!|\!|^\#$.
\end{Theorem}
\begin{proof} By the polar decomposition and Corollary~\ref{C:extension of norm on the completion
algebra}, we may assume that $S$ and $T$ are positive operators. Let
$T_n$ and $S_n$ be two sequences of operators in $\J(\M)$ such that
$\lim_{n\rightarrow\infty}|\!|\!|T-T_n|\!|\!|=\lim_{n\rightarrow\infty}|\!|\!|S-S_n|\!|\!|^\#=0$.
Let $K$ be a positive number such that $|\!|\!|T_n|\!|\!|\leq K$ and
$|\!|\!|S_n|\!|\!|^\#\leq K$ for all $n$ and $\epsilon>0$. Then
there is an $N$ such  that for all $m>n\geq N$,
$|\!|\!|T_m-T_n|\!|\!|<\epsilon/(2K)$ and
$|\!|\!|S_m-S_n|\!|\!|^\#<\epsilon/(2K)$. By Corollary~\ref{C:Holder
inequality for bounded operators},
 $\|T_mS_m-T_nS_n\|_1\leq
\|(T_m-T_n)S_m\|_1+\|T_n(S_m-S_n)\|_1\leq |\!|\!|T_m-T_n|\!|\!|\cdot
|\!|\!|S_m|\!|\!|^\#+|\!|\!|T_n|\!|\!|\cdot
|\!|\!|S_m-S_n|\!|\!|^\#<\epsilon.$ This implies that $\{T_nS_n\}$
is a Cauchy sequence in $\M$ in the $L^1$-norm. Since
$\lim_{n\rightarrow\infty}T_nS_n=TS$ in the measure topology,
$\lim_{n\rightarrow\infty} T_nS_n=TS$ in the $L^1$-norm by
Lemma~\ref{L:convergence in L1 norm}. By Corollary~\ref{C:Holder
inequality for bounded operators}, $\|T_nS_n\|_1\leq
|\!|\!|T_n|\!|\!| \cdot |\!|\!|S_n|\!|\!|^\#$ for every $n$. Hence,
$\|TS\|_1\leq |\!|\!|T|\!|\!| \cdot |\!|\!|S|\!|\!|^\#$.
\end{proof}

Combining Example~\ref{E:Lp space}, \ref{E:Lp for discrete},
\ref{E:dual Lp 1}, \ref{E:dual Lp 2} and Theorem~\ref{T:Holder
inequality}, we obtain the classical non-commutative
H$\ddot{\text{o}}$lder's inequality.

\begin{Corollary}Let $\M$ be an infinite semi-finite factor with the faithful
normal tracial wight $\tau$. If $T\in L^p(\M,\tau)$ and $S\in
L^q(\M,\tau)$, then $TS\in L^1(\M,\tau)$ and
\[\|TS\|_1\leq \|T\|_p\cdot\|S\|_q,
\] where $1\leq p,q\leq \infty$ and $\frac{1}{p}+\frac{1}{q}=1$.
\end{Corollary}

\subsection{Some approximation results}

\begin{Lemma}\label{L:approximation lemma}Let
$|\!|\!|\cdot|\!|\!|$ be a  unitarily invariant norm on $\B(\H)$ and
$T\in \J(\H)$. If $S,S_1,S_2,\cdots, $ are bounded operators in
$\B(\H)$ such that $S=\lim_{n\rightarrow\infty} S_n$ in the strong
operator topology, then
\[\lim_{n\rightarrow\infty} |\!|\!|S_nT-ST|\!|\!|=0.
\]
\end{Lemma}
\begin{proof} Since $T\in \J(\H)$, there is a finite rank projection $E$
such that $T=ETE$. Since $S_n\rightarrow S$ in the strong operator
topology, $S_nE\rightarrow SE$ in the operator topology. By
Proposition~\ref{P:unitarily invariant norms on J(M)},
\[|\!|\!|S_nT-ST|\!|\!|=|\!|\!|S_nET-SET|\!|\!|\leq \|S_nE-SE\|\cdot
|\!|\!|T|\!|\!|\rightarrow 0.
\]
\end{proof}

\begin{Theorem}\label{T:approxmiation theorem} Let
$|\!|\!|\cdot|\!|\!|$ be a  unitarily invariant norm on $\B(\H)$ and
$T\in \overline{\J(\H)_{|\!|\!|\cdot|\!|\!|}}$. If
$S,S_1,S_2,\cdots, $ are bounded operators in $\B(\H)$ such that
$S=\lim_{n\rightarrow\infty} S_n$ in the strong operator topology,
then
\[\lim_{n\rightarrow\infty} |\!|\!|S_nT-ST|\!|\!|=0.
\]
\end{Theorem}
\begin{proof}Since $S_n\rightarrow S$ in the strong operator
topology, there is a number $M>0$ such that $\|S\|\leq M$ and
$\|S_n\|\leq M$ for $n=1,2,\cdots.$  Let $\epsilon>0$. Since $T\in
\overline{\J(\H)_{|\!|\!|\cdot|\!|\!|}}$,  there is a $T'\in \J(\H)$
such that $|\!|\!|T'-T|\!|\!|< \frac{\epsilon}{3M}$. By
Lemma~\ref{L:approximation lemma}, there is an $N$ such that
$|\!|\!|S_nT'-ST'|\!|\!|<\frac{\epsilon}{3}$  for all $n\geq N$. By
Proposition~\ref{P:unitarily invariant norms on J(M)},
\begin{eqnarray}
 \nonumber |\!|\!|S_nT-ST|\!|\!| &\leq & |\!|\!|S_nT-S_nT'|\!|\!|+|\!|\!|S_nT'-ST'|\!|\!|+|\!|\!|ST'-ST|\!|\!| \\
  \nonumber &<& \|S_n\|\cdot|\!|\!|T-T'|\!|\!|+|\!|\!|S_nT'-ST'|\!|\!|+\|S\|\cdot|\!|\!|T'-T|\!|\!| \\
  \nonumber &=&\epsilon.
\end{eqnarray}
\end{proof}

\begin{Corollary}\label{C:approximation corollary} Let
$|\!|\!|\cdot|\!|\!|$ be a  unitarily invariant norm on $\J(\H)$ and
$T\in \overline{\J(\H)_{|\!|\!|\cdot|\!|\!|}}$. If $E_1,E_2,\cdots,
$ are finite rank  projections in $\B(\H)$ such that
$\lim_{n\rightarrow\infty} E_n=I$ in the strong operator topology,
then
\[\lim_{n\rightarrow\infty} |\!|\!|E_nTE_n-T|\!|\!|=0.
\]
\end{Corollary}

Let $\M$ be a diffuse semi-finite von Neumann algebra. Then the set
of finite projections $\{E_\alpha\}$ of $\M$ is a partially ordered
set such that $\lim_{\alpha} E_\alpha=I$ in the strong operator
topology. The following lemma will be useful in the next section.

\begin{Lemma}\label{L:approximation lemma for II infty} Let $\M$ be
a semi-finite von Neumann algebra, and $|\!|\!|\cdot|\!|\!|$ be a
symmetric gauge norm on $\J(\M)$ and $T\in
\overline{\J(\M)_{|\!|\!|\cdot|\!|\!|}}$. Then
\[\lim_{\alpha} |\!|\!|E_\alpha TE_\alpha-T|\!|\!|=0.
\]
\end{Lemma}
\begin{proof} Apply the technique used in the proof of Theorem~\ref{T:approxmiation
theorem}, we need only prove the lemma under the assumption $T\in
\J(\M)$. In this case, the lemma is obvious since there is a
finite projection $E$ such that $ETE=T$.
\end{proof}

\section{Duality and Reflexivity of noncommutative Banach spaces}

 In this section, we assume that $\M$
is an infinite semi-finite factor with a faithful normal tracial
weight $\tau$, $|\!|\!|\cdot|\!|\!|$ is a unitarily invariant norm
on $\J(\M)$ and $|\!|\!|\cdot|\!|\!|^\#$ is the dual unitarily
invariant norm on $\J(\M)$.  We consider the following
two questions in this section:\\

\noindent Question 1: When
$\overline{{\J(\M)}_{|\!|\!|\cdot|\!|\!|^\#}}$ is the dual space of
$\overline{{\J(\M)}_{|\!|\!|\cdot|\!|\!|}}$ in the following sense:
for every $\phi\in \overline{{\J(\M)}_{|\!|\!|\cdot|\!|\!|}}^\#$,
there is a unique $X\in
\overline{{\J(\M)}_{|\!|\!|\cdot|\!|\!|^\#}}$ such that
\[\phi(T)=\tau(TX),\quad\forall\, T\in \overline{{\J(\M)}_{|\!|\!|\cdot|\!|\!|}}
\] and $\|\phi\|=|\!|\!|T|\!|\!|$?\\

\noindent Question 2: When
$\overline{{\J(\M)}_{|\!|\!|\cdot|\!|\!|}}$ is a
reflexive space?\\

 For a projection $E$ in $\M$, we denote by
 $\overline{{(\M_E)}_{|\!|\!|\cdot|\!|\!|}}$ the completion of
 $\M_E=E\M E$ with respect to $|\!|\!|\cdot|\!|\!|$.

\begin{Lemma}\label{L:dual for finite projection} If $\overline{{\J(\M)}_{|\!|\!|\cdot|\!|\!|^\#}}$ is
the dual space of $\overline{{\J(\M)}_{|\!|\!|\cdot|\!|\!|}}$ in the
sense of question 1, then
$\overline{{(\M_E)}_{|\!|\!|\cdot|\!|\!|^\#}}$ is the dual space of
$\overline{{(\M_E)}_{|\!|\!|\cdot|\!|\!|}}$ for every projection $E$
in $\M$ in the sense of question 1.
\end{Lemma}
\begin{proof}Let $\phi\in
\overline{{(\M_E)}_{|\!|\!|\cdot|\!|\!|}}^\#$. We can identify
$\overline{{(\M_E)}_{|\!|\!|\cdot|\!|\!|}}$ with a Banach subspace
of $\overline{{\J(\M)}_{|\!|\!|\cdot|\!|\!|}}$. By the Hahn-Banach
extension theorem, $\phi$ can be extended to a linear functional
$\psi\in \overline{{\J(\M)}_{|\!|\!|\cdot|\!|\!|}}^\#$. By the
assumption of the lemma, there is an operator $X\in
\overline{{\J(\M)}_{|\!|\!|\cdot|\!|\!|^\#}}$ such that for every
operator $T\in \J(\M)$, $\psi(T)=\tau(TX)$.  There is a sequence
of operators $X_n\in \J(\M)$ such that $\lim_{n\rightarrow\infty}
|\!|\!|X_n-X|\!|\!|^\#=0$. Let $Y_n=EX_nE$ and $Y=EXE$. Then
$Y_n\in \M_E$. By Proposition~\ref{P:unitarily invariant implies
trace-preserving},
\[\lim_{n\rightarrow\infty}|\!|\!|Y-Y_n|\!|\!|^\#\leq \lim_{n\rightarrow\infty}
 \|E\|\cdot |\!|\!|X-X_n|\!|\!|^\#\cdot \|E\|=0.
\]
Hence, $Y\in \overline{{(\M_E)}_{|\!|\!|\cdot|\!|\!|^\#}}$. For
every $T\in \M_E$, $\phi(T)=\psi(T)=\tau(TX)=\tau(ETEX)=\tau(TY)$.
This implies that $\overline{{(\M_E)}_{|\!|\!|\cdot|\!|\!|^\#}}$ is
the dual space of $\overline{{(\M_E)}_{|\!|\!|\cdot|\!|\!|}}$ for
every  projection $E$ in $\M$.
\end{proof}

Let $E$ be a (non-zero) finite projection in $\M$. Recall that
$\M_E=E\M E$ is a finite von Neumann algebra with a faithful normal
tracial state $\tau_E$ such that $\tau_E(T)=\frac{\tau(T)}{\tau(E)}$
for $T\in \M_E$. If  $|\!|\!|\cdot|\!|\!|$ is a norm on $\M_E$, the
dual norm of $T\in \M_E$ with respect to $\tau_E$ is defined by
\[|\!|\!|T|\!|\!|^\#_{\M_E,\,\tau_E}=\sup\{|\tau_E(TX)|:\, X\in\M_E,\, |\!|\!|X|\!|\!|\leq 1\}.
\] In the following, we denote by
$\overline{{(\M_E)}_{|\!|\!|\cdot|\!|\!|_{\M_E,\,\tau_E}^\#}}$ the
completion of $\M_E$ with respect to
$|\!|\!|\cdot|\!|\!|_{\M_E,\,\tau_E}^\#$.  By
Lemma~\ref{L:reduction to finite case},
$|\!|\!|T|\!|\!|^\#=\tau(E)|\!|\!|T|\!|\!|_{\M_E,\,\tau_E}^\#$ for
every $T\in \M_E$. Hence,
$\overline{{(\M_E)}_{|\!|\!|\cdot|\!|\!|^\#}}=\overline{{(\M_E)}_{|\!|\!|\cdot|\!|\!|_{\M_E,\,\tau_E}^\#}}$
as linear spaces.\\

The next lemma follows from Theorem~{\bf H} of~\cite{FHNS}.
\begin{Lemma}\label{L:duality for finite vNA} Let $\N$ be a type ${\rm
II}_1$ factor with a faithful normal tracial state $\tau_\N$. Let
$|\!|\!|\cdot|\!|\!|$ be
 a unitarily invariant norm on $\N$ and
$|\!|\!|\cdot|\!|\!|^\#$ be the dual unitarily invariant norm on
$\N$. Let $|\!|\!|\cdot|\!|\!|_1$ be the symmetric gauge norm on
$(L^\infty[0,1],\int_0^1dx)$ corresponding to
$|\!|\!|\cdot|\!|\!|$ on $\N$ as in Theorem~\ref{T:Theorem E}.
Then the following conditions are equivalent:
\begin{enumerate}
 \item $\overline{{\N}_{|\!|\!|\cdot|\!|\!|^\#}}$ is the
dual space of $\overline{{\N}_{|\!|\!|\cdot|\!|\!|}}$ in the sense
of question 1; \item
$\overline{{L^\infty[0,1]}_{|\!|\!|\cdot|\!|\!|_1^\#}}$ is the
dual space of $\overline{{L^\infty[0,1]}_{|\!|\!|\cdot|\!|\!|_1}}$
in the sense of question 1.
\end{enumerate}
\end{Lemma}

The next lemma follows from  Corollary~{\bf 1} of~\cite{FHNS}.

\begin{Lemma}\label{L:Conditional expectation} Let $\N$ be a finite
factor with a faithful normal tracial state $\tau_\N$, and let
$|\!|\!|\cdot|\!|\!|$ be a unitarily invariant norm on $\N$. If
$\A$ is a separable abelian von Neumann subalgebra of $\N$ and
$\mathbf{E}_\A$ is the normal conditional expectation from $\N$
onto $\A$ preserving $\tau_\N$, then
$|\!|\!|\mathbf{E}_{\A}(T)|\!|\!|\leq |\!|\!|T|\!|\!|$ for all
$T\in \N$.
\end{Lemma}

\begin{Definition}\emph{Let $(\M,\tau)$ be a semi-finite von Neumann algebra.
An abelian von Neumann subalgebra $\A$ is called a \emph{normal}
abelian von Neumann subalgebra if $\tau$ is a semi-finite tracial
weight on $\A$}
\end{Definition}

The abelian algebra consists of bounded diagonal operators on a
Hilbert space $\H$ is a normal abelian von Neumann subalgebra of
$\B(\H)$. However, the abelian von Neumann algebra generated by
the bilateral shift operator is not a normal abelian von Neumann
subalgebra of $\B(\H)$.  By~\cite{Ta}, if $\A$ is a normal abelian
von Neumann subalgebra of $\M$, then there is a normal conditional
expectation ${\bf E}_\A$ from $\M$ onto $\A$ preserving $\tau$.\\

The following theorem is the main result of this section.
\begin{Theorem}\label{T:Theorem G}
Let $\M$ be a type ${\rm II}\sb \infty$  factor $($or $\B(\H)$$)$, $|\!|\!|\cdot|\!|\!|$ be
 a unitarily invariant norm on $\M$ and
$|\!|\!|\cdot|\!|\!|^\#$ be the dual unitarily invariant norm on
$\M$. Let $|\!|\!|\cdot|\!|\!|_1$ be the symmetric gauge norm on
$\J(L^\infty[0,\infty))$ $($or $\J(l^\infty(\nn))$ respectively$)$
corresponding to $|\!|\!|\cdot|\!|\!|$ on $\M$ as in
Theorem~\ref{T:Theorem D}. Then the following conditions are
equivalent:
\begin{enumerate}
 \item $\overline{{\J(\M)}_{|\!|\!|\cdot|\!|\!|^\#}}$ is the
dual space of $\overline{{\J(\M)}_{|\!|\!|\cdot|\!|\!|}}$ in the
following sense:  $\forall\,\psi\in
\overline{{\J(\M)}_{|\!|\!|\cdot|\!|\!|}}^\#$, there is a unique
$X\in \overline{{\J(\M)}_{|\!|\!|\cdot|\!|\!|}^\#}$ such that
$\psi(T)=\tau(TX)$ for all $T\in
\overline{{\J(\M)}_{|\!|\!|\cdot|\!|\!|}}$ and
$\|\psi\|=|\!|\!|X|\!|\!|^\#$; \item
$\overline{{\J(L^\infty[0,\infty))}_{|\!|\!|\cdot|\!|\!|_1^\#}}$
 $($or $\overline{{\J(l^\infty(\nn))}_{|\!|\!|\cdot|\!|\!|_1^\#}}$ respectively$)$ is the
dual space of
$\overline{{\J(L^\infty[0,\infty))}_{|\!|\!|\cdot|\!|\!|_1}}$
 $($or $\overline{{\J(l^\infty(\nn))}_{|\!|\!|\cdot|\!|\!|_1}}$ respectively$)$in the
same sense.
\end{enumerate}
\end{Theorem}
\begin{proof} We prove the theorem for type
${\rm II}\sb \infty$ factors. The proof for type ${\rm I}\sb \infty$
factors
is similar.\\

\, $``1\Rightarrow 2"$.  By Theorem~\ref{T:Theorem D} and
Lemma~\ref{L:isomorphism}, there is a normal separable diffuse
abelian von Neumann subalgebra $\A$ of $\M$ and a
$\ast$-isomorphism $\alpha$ from $\A$ onto $L^\infty[0,\infty)$
such that $\tau=\int_0^\infty\,dx\cdot \alpha$ and
$|\!|\!|\alpha(T)|\!|\!|_1=|\!|\!|T|\!|\!|$ for every $T\in
\J(\A)$.  a By Theorem~\ref{T:Theorem E},
$|\!|\!|\alpha(T)|\!|\!|_1^\#=|\!|\!|T|\!|\!|^\#$ for every $T\in
\J(\A)$. So we need only  prove that
$\overline{{\J(\A)}_{|\!|\!|\cdot|\!|\!|^\#}}$ is the dual space
of $\overline{{\J(\A)}_{|\!|\!|\cdot|\!|\!|}}$ in the sense of
question 1. Let $\phi\in
\overline{{\J(\A)}_{|\!|\!|\cdot|\!|\!|}}^\#$. By the Hahn-Banach
extension theorem, $\phi$ can be extended to a bounded linear
functional $\psi$ on $\overline{{\J(\M)}_{|\!|\!|\cdot|\!|\!|}}$
such that $\|\psi\|=\|\phi\|$. By the assumption of the theorem,
there is an operator $X\in
\overline{{\J(\M)}_{|\!|\!|\cdot|\!|\!|}^\#}$ such that
$\psi(S)=\tau(SX)$ for all $S\in
\overline{{\J(\M)}_{|\!|\!|\cdot|\!|\!|}}$ and
$\|\psi\|=|\!|\!|X|\!|\!|^\#$.  There is a sequence of operators
$X_n\in \J(\M)$ such that
$\lim_{n\rightarrow\infty}|\!|\!|X-X_n|\!|\!|=0$. Let ${\bf E}_\A$
be the normal conditional expectation from $\M$ onto $\A$ such
that $\tau(T)=\tau({\bf E}_\A(T))$ for all $T\in \J(\M)$.  Let
$Y_n=\mathbf{E}_\A(X_n)$. By Lemma~\ref{L:Conditional
expectation}, $|\!|\!|Y_n|\!|\!|\leq |\!|\!|X_n|\!|\!|$.
Therefore,
 $\{Y_n\}$
is a Cauchy sequence in $\J(\A)$ with respect to norm
$|\!|\!|\cdot|\!|\!|$. Let $Y=\lim_{n\rightarrow\infty} Y_n$ in
$\overline{{\J(\A)}_{|\!|\!|\cdot|\!|\!|}}$. Then for every $S\in
\J(\A)$,
\[\phi(S)=\psi(S)=\tau(SX)=\lim_{n\rightarrow\infty}\tau(SX_n)=
\lim_{n\rightarrow\infty}\tau(\mathbf{E}_\A
(SX_n))=\lim_{n\rightarrow\infty}\tau(SY_n)=\tau(SY).
\] By Lemma~\ref{L:dual norms for unbounded operators},
$\|\phi\|=|\!|\!|Y|\!|\!|_\A^\#$.\\

$``2\Rightarrow 1"$. Let $\psi\in
\overline{{\J(\M)}_{|\!|\!|\cdot|\!|\!|}^\#}$, and let $\A$ be a
normal separable diffuse abelian von Neumann subalgebra of $\M$
and $\phi$ be the restriction of $\psi$ to $\J(\A)$. By
Lemma~\ref{L:isomorphism}, we can identify $(\A,\tau)$ with
$(L^\infty[0,\infty), \int_0^\infty\,dx)$. By the assumption of
the theorem,
$\overline{{\J(L^\infty[0,\infty))}_{|\!|\!|\cdot|\!|\!|_1^\#}}$
is the dual space of
$\overline{{\J(L^\infty[0,\infty))}_{|\!|\!|\cdot|\!|\!|_1}}$.
Hence, $\overline{(\A)_{|\!|\!|\cdot|\!|\!|^\#}}$  is the dual
space of $\overline{(\A)_{|\!|\!|\cdot|\!|\!|}}$. Therefore, there
is an operator $Y\in \overline{(\A)_{|\!|\!|\cdot|\!|\!|^\#}}$
such that $\phi(T)=\tau(TY)$ for all $T\in \J(\A)$ and
$|\!|\!|\phi|\!|\!|^\#=|\!|\!|Y|\!|\!|^\#$.\\

Let $E\in \A$ be a finite projection. By Lemma~\ref{L:dual for
finite projection}, $\overline{(\A_E)_{|\!|\!|\cdot|\!|\!|_1^\#}}$
is the dual space of $\overline{(\A_E)_{|\!|\!|\cdot|\!|\!|_1}}$.
By Lemma~\ref{L:reduction to finite case} and Lemma~\ref{L:duality
for finite vNA}, $\overline{(\M_E)_{|\!|\!|\cdot|\!|\!|^\#}}$  is
the dual space of $\overline{(\M_E)_{|\!|\!|\cdot|\!|\!|}}$. So
there is a unique operator $X_E\in
\overline{(\M_E)_{|\!|\!|\cdot|\!|\!|^\#}}$ such that for all
$T\in \M_E$, $\psi(T)=\tau(TX_E)$.  Define $\psi_E(T)=\psi(T)$ for
$T\in \overline{{\M_E}_{|\!|\!|\cdot|\!|\!|}}$ and
$\phi_E(S)=\phi(S)$ for $S\in
\overline{{\A_E}_{|\!|\!|\cdot|\!|\!|}^\#}$. Then
$\phi_E(S)=\phi(S)=\tau(SY)=\tau(ESEY)=\tau(SEYE)$. By
Lemma~\ref{L:duality for finite vNA},
$|\!|\!|X_E|\!|\!|^\#=\|\psi_E\|=\|\phi_E\|=|\!|\!|EYE|\!|\!|^\#$.\\

Let $E,F$ be two finite projections in $\A$.  Then  $
|\!|\!|X_E-X_F|\!|\!|^\#=\|\psi_E-\psi_F\|=\|\phi_E-\phi_F\|=|\!|\!|EYE-FYF|\!|\!|^\#$.
By Lemma~\ref{L:approximation lemma for II infty}, $\lim_\alpha
|\!|\!|E_\alpha YE_\alpha-Y|\!|\!|^\#=0$.
 So $\{X_{E_\alpha}\}$ is
a Cauchy sequence in $\overline{{\J(\M)}_{|\!|\!|\cdot|\!|\!|}^\#}$.
Let $X=\lim_\alpha {E_\alpha}XE_\alpha$. Then $X\in
\overline{{\J(\M)}_{|\!|\!|\cdot|\!|\!|}^\#}$ and $\psi(T)=\tau(TX)$
for all $T\in \M_E$ and finite projections $E\in \A$. \\

Let  $\B$ be another normal separable diffuse abelian von Neumann
subalgebra of $\M$. Similar arguments as above shows that there is
an operator $Z\in \overline{{\J(\M)}_{|\!|\!|\cdot|\!|\!|}^\#}$
such that $\psi(T)=\tau(TZ)$ for all $T\in \M_{E'}$ and finite
projections $E'\in \B$. Let $E\in \A$ and $E'\in\B$ be finite
projections. Then $F=E\wedge E'$ is a finite projection and
$\psi(T)=\tau(TFXF)=\tau(TFZF)$ for all $T\in \M_F$. Therefore,
$FXF=FZF$. Since we can choose $E_n\in\A$ and $E_n'\in \B$ such
that $F_n=E_n\wedge E_n'\rightarrow 1$ in the strong operator
topology, $F_nXF_n=F_nZF_n$ implies that $X=Z$ by Theorem 3
of~\cite{Ne}.\\

Note that if  $T\in \J(\M)$ is a positive operator, then there is
a normal separable diffuse abelian von Neumann subalgebra of $\M$
that contains $T$. Therefore, $\psi(T)=\tau(TX)$ for all positive
operators $T\in \J(\M)$ and hence for all operators in $\J(\M)$.
Since $X\in \overline{{\J(\M)}_{|\!|\!|\cdot|\!|\!|}^\#}$, by
Theorem~\ref{T:Holder inequality}, $\psi(T)=\tau(TX)$ for all
$T\in \overline{{\J(\M)}_{|\!|\!|\cdot|\!|\!|}}$. By
Lemma~\ref{L:dual norms for unbounded operators},
$\|\psi\|=|\!|\!|X|\!|\!|^\#$.
\end{proof}

\begin{Corollary} For $1\leq p<\infty$, $L^q(\M,\tau)$ is the dual
space of $L^q(\M,\tau)$ in the sense of Question 1, where
$\frac{1}{p}+\frac{1}{q}=1$.
\end{Corollary}

The following theorem is a corollary of Theorem~\ref{T:Theorem G}.
\begin{Theorem}\label{T:Theorem H} Let $\M$ be a type ${\rm II}\sb
\infty$ factor $($or $\B(\H)$$)$, $|\!|\!|\cdot|\!|\!|$ be
 a unitarily invariant norm on $\J(\M)$  and $|\!|\!|\cdot|\!|\!|_1$ be the symmetric gauge norm on
$\J(L^\infty[0,\infty))$ $($or $\J(l^\infty(\nn))$ respectively$)$
corresponding to $|\!|\!|\cdot|\!|\!|$ on $\J(\M)$ as in
Theorem~\ref{T:Theorem D}. Then the following conditions are
equivalent:
\begin{enumerate}
\item $\overline{{\J(\M)}_{|\!|\!|\cdot|\!|\!|}}$ is a reflexive
space; \item
$\overline{{\J(L^\infty[0,\infty))}_{|\!|\!|\cdot|\!|\!|_1}}$
$($or $\J(l^\infty(\nn))$ respectively$)$ is a reflexive space.
\end{enumerate}
\end{Theorem}

\begin{Corollary} For $1<p<\infty$, $L^p(\M,\tau)$ is a reflexive space.
\end{Corollary}

\begin{Example}\emph{By Theorem~\ref{T:dual norms of Ky Fan norms} and Theorem~\ref{T:Theorem H},
 for $0\leq t\leq \infty$,
$\overline{{\M}_{|\!|\!|\cdot|\!|\!|_{(t)}}}$ is not a reflexive
space.}
\end{Example}

\vspace{.2in} \noindent Department of Mathematics, University of New
Hampshire

\noindent Durham, NH  03824

\noindent {\em E-mail address: } [Junsheng Fang]\,\,
jfang\@@cisunix. unh.edu

\noindent {\em E-mail address: } [Don Hadwin]\quad\quad
don\@@math.unh.edu

\end{document}